\documentclass[11pt]{amsart}

\usepackage[a4paper,hmargin=3cm,vmargin=3cm]{geometry}
\usepackage{amsfonts,amssymb,amscd,amstext}
\usepackage{graphicx}
\usepackage[dvips]{epsfig}
\usepackage[latin1]{inputenc}

\usepackage{fancyhdr}
\usepackage{times}

\usepackage{enumerate}
\usepackage{titlesec}
\usepackage{mathrsfs}
\usepackage{stmaryrd}

\def\R{\mathbb{R}}

\def\N{\mathbb{N}}

\def\C{\mathbb{C}}

\newcommand{\ben}{\begin{enumerate}}
\newcommand{\bit}{\begin{itemize}}
\newcommand{\een}{\end{enumerate}}
\newcommand{\eit}{\end{itemize}}

\newcommand{\ed}{\end{document}}

\def\cA{\mathcal{A}}
\def\cU{\mathcal{U}}

\def\cS{\mathcal{S}}
\def\cC{\mathcal{C}}
\def\cD{\mathcal{D}}

\def\cW{\mathcal{W}}
\def\cS{\mathcal{S}}
\def\cV{\mathcal{V}}

\def\cH{\mathcal{H}}
\def\cL{\mathcal{L}}
\def\cM{\mathcal{M}}

\def\cG{\mathcal{G}}

\def\cN{\mathcal{N}}

\def\cK{\mathcal{K}}
\def\cO{\mathcal{O}}

\let\8=\infty \let\0=\emptyset

\let\landa=\lambda
\let\alfa=\alpha

\let\parc=\partial

\def\ep{\varepsilon}

\def\landa{\lambda}

\def\flecha{\rightarrow}
\def\esiz{\langle}
\def\esde{\rangle}

\def\S{\Sigma}

\def\cte.{\mathop{\rm cte.}\nolimits}

\def\N{\mathbb{N}}

\def\R{\mathbb{R}}

\def\C{\mathbb{C}}

\def\S{\mathbb{S}}

\newfont{\bb}{msbm10 at 12pt}

\titleformat{\section}
{\filcenter\bfseries\large} {\thesection{.}}{0.2cm}{}
\titleformat{\subsection}[runin]
{\bfseries} {\thesubsection{.}}{0.15cm}{}[.]
\titleformat{\subsubsection}[runin]
{\em}{\thesubsubsection{.}}{0.15cm}{}[.]

\usepackage[up,bf]{caption}

\newtheorem{theorem}{Theorem}[section]

\newtheorem{remark}[theorem]{Remark}

\newtheorem{definition}[theorem]{Definition}

\newtheorem{assertion}[theorem]{Assertion}

\theoremstyle{definition}

\numberwithin{equation}{section}
\numberwithin{figure}{section}

\usepackage{color}

\begin{document}
\fancyhead[LO]{Homogeneous solutions of degenerate elliptic equations}
\fancyhead[RE]{José A. Gálvez, Pablo Mira}
\fancyhead[RO,LE]{\thepage}

\thispagestyle{empty}

\begin{center}
{\bf \LARGE Linearity of homogeneous solutions to degenerate\\[0.2cm] elliptic equations in dimension three}
\vspace*{5mm}

\hspace{0.2cm} {\Large José A. Gálvez and Pablo Mira}
\end{center}

\footnote[0]{
\mbox{} \hspace{0.25cm} \emph{Keywords}: Degenerate elliptic equation, homogeneous function, Alexandrov conjecture, saddle  function.}



\vspace*{7mm}

\begin{quote}
{\small
\noindent {\bf Abstract}\hspace*{0.1cm}
Given a linear elliptic equation $\sum a_{ij} u_{ij} =0$ in $\R^3$, it is a classical problem to determine if its degree-one homogeneous solutions $u$ are linear. The answer is negative in general, by a construction of Martinez-Maure. In contrast, the answer is affirmative in the uniformly elliptic case, by a theorem of Han, Nadirashvili and Yuan, and it is a known open problem to determine the degenerate ellipticity condition on $(a_{ij})$ under which this theorem still holds. In this paper we solve this problem. We prove the linearity of $u$ under the following degenerate ellipticity condition for $(a_{ij})$, which is sharp by Martinez-Maure example: if $\cK$ denotes the ratio between the largest and smallest eigenvalues of $(a_{ij})$, we assume $\cK|_{\cO}$ lies in $L_{\rm loc}^1$ for some connected open set $\cO\subset \S^2$ that intersects any configuration of four disjoint closed geodesic arcs of length $\pi$ in $\S^2$. Our results also give the sharpest possible version under which an old conjecture by Alexandrov, Koutroufiotis and Nirenberg (disproved by Martinez-Maure's example) holds.

\vspace*{0.1cm}

}
\end{quote}


\section{Introduction}

Let $u\in C^2(\R^3\setminus\{0\})$ be a degree-one (positively) homogeneous solution to the linear equation
\begin{equation}\label{edp1}
\sum_{i,j=1}^3 a_{ij} u_{ij} =0, \hspace{1cm} a_{ij}\in L^{\8}(\R^3),
\end{equation}
 in $\R^3$, i.e., $u(\rho x)=\rho u(x)$ for all $\rho>0$, $x\in \R^3$. Assume that \eqref{edp1} is elliptic, i.e., 
\begin{equation}\label{elipp}
\text{$(a_{ij}(x))$ is positive definite}
\end{equation} for every $x\in \R^3$. Note that the $a_{ij}$ are not continuous. \emph{Must then $u$ be a linear function?}

\begin{figure}[htbp]
    \centering
    \includegraphics[width=0.9\textwidth]{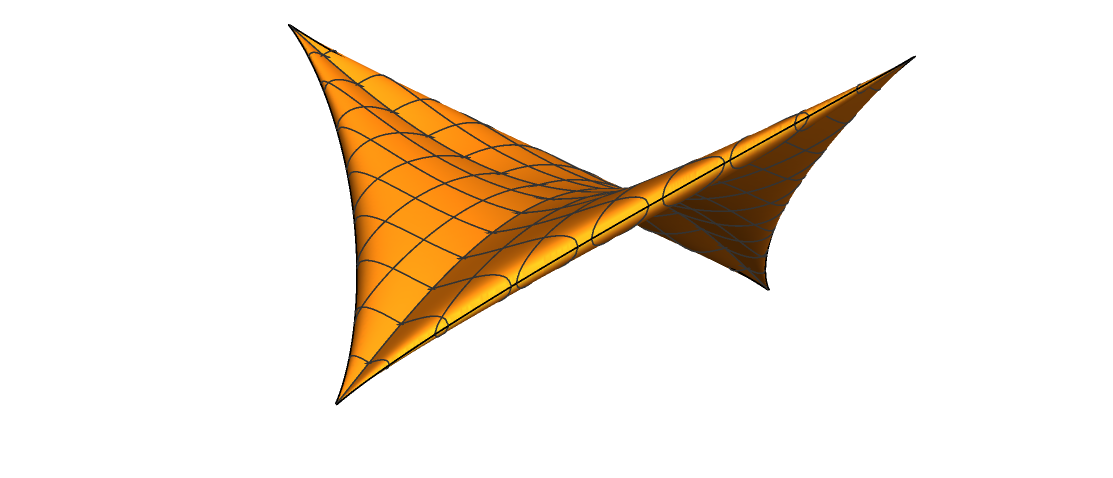}
     \caption{Martinez-Maure's hegdehog $\nabla u(\S^2)$, where $u$ solves \eqref{edp1}-\eqref{elipp}. The preimage in $\S^2$ of each of the four \emph{horns} of the example is a geodesic semicircle.} 
\label{fig:crosscap}
\end{figure}

This is a classical question motivated by global surface theory. Using an equivalent formulation, Alexandrov proved in 1939 that the answer is affirmative if $u$ is real analytic (\cite{A0}), and conjectured that an affirmative answer should also hold in the general case (\cite{A}, p. 352). The validity of this conjecture remained elusive for a long time, until Martinez-Maure \cite{MM} constructed in 2001 a striking  $C^2$ counterexample to it. Specifically, he proved the existence of a nonlinear function $h\in C^2(\S^2)$ such that the \emph{hedgehog} $\psi (\nu):= \nabla h (\nu) + h(\nu) \nu : \S^2\flecha \R^3$ has negative curvature at its regular points. The homogeneous extension of degree one $u$ to $\R^3$ of $h$ gives a counterexample to Alexandrov's conjecture. See Figure \ref{fig:crosscap}.

In contrast, in 2003 Han, Nadirashvili and Yuan \cite{HNY} proved that the Alexandrov conjecture holds in the uniformly elliptic case. This solved an open problem by Safonov \cite{Sa}. Specifically, if $0< \landa(x)\leq \Lambda(x)$ are the smallest and largest eigenvalues of $(a_{ij}(x))$, and we denote $\cK(x):= \Lambda(x)/\landa(x)\geq 1$, Han, Nadirashvili and Yuan imposed the condition
\begin{equation}\label{unife}
 \cK\in L^{\8}(\R^3) 
\end{equation}
and proved the following remarkable result:
\begin{theorem}[\cite{HNY}]\label{th:nadi}
Any $1$-homogeneous solution $u\in W^{2,2}_{\rm loc} (\R^3)$ to \eqref{edp1}-\eqref{unife}
is linear.
\end{theorem}
An alternative proof of Theorem \ref{th:nadi} was obtained in 2016 by Guan, Wang and Zhang \cite{GWZ}, again under very weak regularity assumptions on $u$. For that, they treated the problem directly as a uniformly elliptic equation in $\S^2$, and gave an elegant argument using the Bers-Nirenberg unique continuation theorem. A different approach to Theorem \ref{th:nadi} via Poincaré-Hopf index theory was given by the authors and Tassi in \cite{GMT}. The problem of the linearity of homogeneous solutions to \eqref{edp1}-\eqref{elipp} is discussed in detail in the book \cite{NTV} by Nadirashvili, Tkachev and Vladut. 

The uniform ellipticity assumption \eqref{unife} in Theorem \ref{th:nadi} cannot be weakened to \emph{plain} ellipticity \eqref{elipp}, by Martinez-Maure's example. A known natural open problem proposed by Guan, Wang and Zhang (see \cite[Remark 8]{GWZ}) is to establish what degenerate ellipticity conditions on the coefficients $a_{ij}$ are sufficient for Theorem \ref{th:nadi} to hold, even when $u$ is smooth. 

In this paper we give an answer to this problem. We explain next our main results.

Let $u\in C^2(\R^3\setminus\{0\})$ be a degree-one homogeneous solution to a linear equation \eqref{edp1}. 
By homogeneity, $u$ also satisfies \eqref{edp1} for the coefficients $a_{ij}= a_{ij}(x/|x|)$. For this reason, our hypotheses on $(a_{ij})$ will be directly viewed at points $x\in \S^2$. Instead of \eqref{unife}, we will just assume that the considerably weaker condition
\begin{equation}\label{coe}
\cK |_{\cO}\in L_{\rm loc}^1(\cO)
\end{equation}
holds for \emph{some} connected open set $\cO\subset \S^2$ that it intersects any configuration of four disjoint \emph{geodesic semicircles} (i.e. closed geodesic arcs of length $\pi$) in $\S^2$. Such a set $\cO$ can be quite small. For instance, $\cO$ can be chosen as any connected open set of $\S^2$ that contains an arbitrarily thin collar along a geodesic, $C_{\gamma}:=\{x\in \S^2 : \esiz x,\nu_0\esde \in (0,\ep)\}$ for some $\nu_0\in \S^2$, $\ep>0$. 

We prove:

\begin{theorem}\label{th:maine}
Any $1$-homogeneous solution $u\in C^2(\R^3\setminus\{0\})$ to \eqref{edp1}, \eqref{elipp}, \eqref{coe} is linear. 
\end{theorem}

The four semicircles condition imposed on $\cO$ is sharp.
Indeed, Martinez-Maure's example in \cite{MM} yields a $1$-homogeneous function $u\in C^2(\R^3\setminus\{0\})$ such that $D^2 u$ is indefinite whenever it is non-zero, and so that $\{x \in \S^2 : D^2 u (x)=0\}$ agrees exactly with a certain configuration $\Gamma\subset \S^2$ of four disjoint geodesic semicircles. By the indefinite nature of $D^2 u$, we can view $u$ as a solution to some elliptic equation \eqref{edp1}-\eqref{elipp}, and the related function $\cK$ associated to the coefficients $(a_{ij})$ of this equation lies in $L_{\rm loc}^1(\cO)$ for any open set $\cO\subset \S^2$ disjoint from $\Gamma$. 

We can actually prove a more general version of Theorem \ref{th:maine}, that holds under degenerate ellipticity conditions: 

\begin{theorem}\label{th:main}
Let $u\in C^2(\R^3\setminus\{0\})$ be a $1$-homogeneous solution to \eqref{edp1}. Assume 
\begin{equation}\label{co}
\left\{\def\arraystretch{1.5}\begin{array}{l}
\text{(i) $(a_{ij}(x))$ is semi-positive definite $\forall x\in \R^3$.}\\ \text{(ii) The restriction of $(a_{ij}(x))$ to the plane $x^{\perp}$ is non-zero, $\forall x\in \R^3-\{0\}$.}\\ \text{(iii) $(a_{ij})$ is positive definite a.e. on $\cO$, and $\cK |_{\cO}\in L_{\rm loc}^1(\cO)$},\end{array} \right.
\end{equation}where, in $(iii)$, $\cO\subset \S^2$ is some connected open set that intersects any configuration of four disjoint geodesic semicircles.
Then, $u$ is linear.
\end{theorem}
Note that (i) extends \eqref{elipp} to the degenerate elliptic setting, and (ii) is needed in that general context to ensure that  \eqref{edp1} is non-trivial when restricted to $1$-homogeneous functions.

The proof of Theorem \ref{th:main} is a blend of geometric and analytic arguments, and is presented in Section \ref{sec:2}. The idea, following Alexandrov \cite{A0}, is to show that $\nabla u (\S^2)$ reduces to a point, by analyzing the support planes in $\R^3$ of this compact set. In the uniformly elliptic case, Han, Nadirashvili and Yuan \cite{HNY} used this idea and the maximum principle to show that $\nabla u(\S^2)$ is a point. In our situation given by \eqref{co}, we will use instead the Stoilow factorization for planar mappings of finite distortion \cite{IS,AIM}. However, the main difficulty of the proof is that we are not assuming that $\cK\in L^1 (\S^2)$, but only that its restriction to the possibly quite small set $\cO\subset \S^2$ lies in $L_{\rm loc}^1$. In order to deal with this general situation, we will use an idea of Pogorelov \cite{P1}. In \cite{P1}, Pogorelov claimed a proof of Alexandrov's conjecture, something that is incorrect by the example in \cite{MM}. Pogorelov's argument was based on the deep idea of controlling the connected components in which some suitable planes of $\R^3$ divide the saddle graph $\Sigma$ in $\R^3$ given by $z=u(x,y,1)$. However, this is a delicate question, and the short argument presented in \cite{P1} has several errors in the way these connected components are handled (one of them was pointed out in \cite{Pan1}). Our proof of Theorem \ref{th:main} springs from Pogorelov's brilliant idea, but we give a different, subtler argument that yields full control of the connected components mentioned above.

The term \emph{Alexandrov conjecture} is often used in the literature in reference to a more general statement, in which \eqref{edp1} is allowed to be degenerate elliptic; see e.g. \cite{MM,Pan1,NTV}. This conjecture admits several equivalent formulations, one of which is the following one, proposed in 1973 by Koutroufiotis and Nirenberg \cite{K}:

{\bf The Alexandrov-Koutroufiotis-Nirenberg conjecture:} \emph{Any $C^2$ function $v$ in $\S^2$ that satisfies ${\rm det}(\nabla_{\S^2}^2 v)\leq 0$ at every point must be linear, i.e., $\nabla_{\S^2}^2 v=0$.}

Here, as usual, the spherical Hessian $\nabla_{\S^2}^2v $ is defined by $\nabla_{\S^2}^2v(q)=(v_{ij}(q) + v(q)\delta_{ij})$, where $v_{ij}$ are covariant derivatives with respect to a local orthonormal frame in $\S^2$, see e.g. \cite{GWZ}. We say that $v\in C^2(\S^2)$ is a \emph{saddle function} on $\S^2$ if it satisfies ${\rm det}(\nabla_{\S^2}^2 v)\leq 0$. The conjecture is then that saddle functions on $\S^2$ are linear.

The support function $h$ of Martinez-Maure's hedgehog in \cite{MM} gives a $C^2$ counterexample to this conjecture. Panina's construction in \cite{Pan1} provides $C^{\8}$ counterexamples, which are actually linear in large open regions of $\S^2$. Based on these results, Nadirashvili, Tkachev and Vladut proposed in \cite[Conjecture 1.6.1]{NTV} a \emph{lopped version} of the conjecture, which can be rephrased as follows: \emph{any $C^2$ saddle function on $\S^2$ is linear in some open set}.

This beautiful conjecture in \cite{NTV} is open if $v$ is at least of class $C^3$, but in the general $C^2$ category, one should reformulate it slightly. Indeed, Martinez-Maure's saddle function $h\in C^2(\S^2)$ satisfies that $\{q\in \S^2 : \nabla_{\S^2}^2 h(q)=0\}$ is the union of four disjoint geodesic semicircles; in particular, $h$ is not linear on any open set of $\S^2$. Thus, the best possible \emph{lopped} conjecture that can hold in the general $C^2$ case is that any saddle function $v\in C^2(\S^2)$ always satisfies $\nabla_{\S^2}^2 v=0$ along four disjoint geodesic semicircles. We will prove this exact result as a part of our proof of Theorem \ref{th:main}; see Section \ref{sec:4s}.
\begin{theorem}\label{th:4sh}
Let $v\in C^2(\S^2)$ satisfy ${\rm det}(\nabla_{\S^2}^2 v)\leq 0$. Then $\nabla_{\S^2}^2 v =0$ along four disjoint geodesic semicircles of $\S^2$.
\end{theorem}
Theorem \ref{th:4sh} gives then the sharpest possible version for which the conjecture by Alexandrov, Koutroufiotis and Nirenberg is true, i.e., the sharpest possible \emph{linearity} theorem for saddle $C^2$ functions in $\S^2$. We should note that Panina claimed in \cite{Pan2} a very general statement that would have Theorem \ref{th:4sh} as a particular case. However, the very short argument given in \cite{Pan2} is not correct; for instance, it relies on Pogorelov's incorrect study of the connected components problem. In Theorem \ref{th:4s} we will give an alternative formulation of Theorem \ref{th:4sh}, in the context of the Weingarten inequality $(\kappa_1-c)(\kappa_2-c)\leq 0$ for ovaloids of $\R^3$.

The Alexandrov conjecture has been linked by Mooney \cite{Mo} to the existence of Lipschitz minimizers to functionals $\int F(\nabla u) dx$ in $\R^3$, with $F$ strictly convex, that are $C^1$ except at a finite number of points. It has also been linked in \cite{HNY,NTV,NV} to the classification of degree-two homogeneous solutions to elliptic Hessian equations $F(D^2 u)=0$ in $\R^3$. In particular, our results here might be of interest regarding the following conjecture in the book by Nadirashvili, Tkachev and Vladut, see \cite[Conjecture 1.6.3]{NTV}: \emph{a degree-two homogeneous smooth solution $u$ to a degenerate elliptic Hessian equation $F(D^2 u)=0$ in $\R^3$ must be a quadratic polynomial}.

The authors are grateful to Yves Martinez-Maure for enlightening comments and discussions.
\section{Proof of Theorem \ref{th:main}}\label{sec:2}

Let $u\in C^2(\R^3\setminus\{0\})$ be a degree-one homogeneous solution to \eqref{edp1}, where \eqref{co} holds. We will assume throughout the proof that $u$ is not linear, i.e. $D^2 u$ is not identically zero on $\R^3$, and reach a contradiction. We will split the proof into several steps.

\vspace{0.2cm}

{\bf Step 1:} \emph{Connection with quasiregular mappings}. 

In this step we relate the conditions in \eqref{co} with the theory of planar mappings with finite distortion, in order to apply the Stoilow factorization by Iwaniec-Sverak \cite{IS} to our context.

\vspace{0.2cm}

Consider arbitrary Euclidean coordinates $(x,y,z)$ in $\R^3$ centered at the origin, and define $h\in C^2(\R^2)$ by
\begin{equation}\label{eq:0}
h(x,y):= u(x,y,1).
\end{equation}
Note that $u(x,y,z)=z h(x/z,y/z)$ for all $z>0$, by homogeneity.  Then we have (see \cite{HNY})
\begin{equation}\label{eq:1g}
\nabla u (x,y,1)= (h_x,h_y, h -x h_x-y h_y)
\end{equation}
and 
\begin{equation}\label{eq:1h}
D^2 u (x,y,1)= \left(\def\arraystrecth{1.5}\begin{array}{ccl} h_{xx} & h_{xy} & -x h_{xx} - y h_{xy} \\ * & h_{yy} & -x h_{xy} - y h_{yy} \\ * & * & x^2 h_{xx} + 2 x y h_{xy} + y^2 h_{yy} \end{array} \right).
\end{equation}
From here and the invariance of \eqref{edp1} by Euclidean isometries we see that the restriction of \eqref{edp1} to points of the form $(x,y,1)$ turns into a linear PDE for $h$,
\begin{equation}\label{eq:2}
A_{11}h_{xx}+2 A_{12}h_{xy}+A_{22}h_{yy} =0.
\end{equation}
Specifically, if we denote $\cA:=(a_{ij}(x,y,1))$ and $\cM:= (A_{ij}(x,y))$, by \eqref{eq:1h}, the coefficients of \eqref{eq:2} are given for $i,j\in \{1,2\}$ by
\begin{equation}\label{eq:3}
A_{ij} = w_i \cdot \cA \cdot w_j^T,
\end{equation}
where $w_1:=(1,0,-x)$ and $w_2:=(0,1,-y)$. In other words, the bilinear form defined by $\cM$ is the restriction of the one given by $\cA$ to the plane of $\R^3$ orthogonal to $(x,y,1)$. By (i) and (ii) in \eqref{co}, the matrix $\cM$ is semi-positive definite and non-zero for all $(x,y)$. This clearly implies by \eqref{eq:2} that, for any $(x,y)$, 
\begin{equation}\label{eq:11}
h_{xx} h_{yy}-h_{xy}^2\leq 0 .
\end{equation}
The converse of this property also holds, i.e., if $h(x,y)$ satisfies \eqref{eq:11}, it solves a degenerate elliptic equation \eqref{eq:2} in $\R^2$, for adequate coefficients $A_{ij}$; see e.g. \cite{Pu} for a similar argument in the elliptic case. Hence, if for \emph{any} Euclidean linear coordinate system $(x,y,z)$, the function $h(x,y)$ given by \eqref{eq:0} satisfies \eqref{eq:11}, then $u$ solves a linear equation \eqref{edp1} whose coefficients $a_{ij}$ satisfy (i), (ii) in \eqref{co}. 
 
Consider the smallest and largest eigenvalues $\landa\leq \Lambda$ among the three eigenvalues of $\cA$ at $(x,y,1)$, and let $\landa_1\leq \landa_2$ denote the eigenvalues of $\cM$. By \eqref{eq:3}, we have 
 \begin{equation}\label{eq:4}
 0\leq\landa \leq \landa_1\leq \landa_2 \leq \Lambda.
 \end{equation}

Choose next a point $\nu_0\in \cO\subset\S^2$ with positive $z$-coordinate, and express it as
 \begin{equation}\label{eq:9a}
\nu_0 = \frac{1}{\sqrt{1+x_0^2+y_0^2}}(x_0,y_0,1).
\end{equation} 
Since $(a_{ij})$ is positive definite a.e. on $\cO$ by (iii) in \eqref{co}, the matrix $\cM$ is positive definite a.e. around $(x_0,y_0)$, by \eqref{eq:4}. Dividing by $A_{11}+A_{22}$, we can rewrite \eqref{eq:2} as
 \begin{equation}\label{eq:5}
2 h_{w\bar{w}} + \mu h_{ww} + \overline{\mu} h_{\bar{w}\bar{w}} =0
 \end{equation}
around $w_0:=x_0+iy_0$, where $w=x+iy$ and 
\begin{equation}\label{eq:6}
\mu = \frac{A_{11}-A_{22}+2i A_{12}}{A_{11} +A_{22}}.
\end{equation}
Thus,
\begin{equation}\label{eq:7}
|\mu|=\frac{K_{\mu}-1}{K_{\mu}+1}<1,\hspace{0.5cm} \text{ where } \hspace{0.5cm} K_{\mu}:= \frac{\landa_2}{\landa_1}\geq 1.
\end{equation}
If we now write $f:=h_w$, then by \eqref{eq:5} and \eqref{eq:7} we have 
 \begin{equation}\label{eq:8}
 |f_{\bar{w}}|\leq |\mu| |f_w|, \hspace{1cm} |\mu|<1 \text{ a.e. around $w_0$}.
 \end{equation}
Let us control next the dilatation quotient of $f$. If we denote $$J(w,f):=|f_w|^2 - |f_{\bar{w}}|^2 \geq 0, \hspace{1cm} |Df(w)|:= |f_w|+ |f_{\bar{w}}|,$$ the dilatation quotient of $f$ is given for any $w\in \C$ with $J(w,f)\neq 0$ by $$K(w,f)=\frac{|Df(w)|^2}{J(w,f)}\geq 1.$$ At the points where $|Df(w)|=J(w,f)=0$, we define $K(w,f):=1$. Thus, $K(w,f)$ is defined a.e. around $w_0$, and by \eqref{eq:7} and \eqref{eq:8} we have at points with $J(w,f)\neq 0$
\begin{equation}\label{eq:9}
K(w,f)\leq \frac{(|f_w|+|\mu| |f_w|)^2}{|f_w|^2-|\mu|^2 |f_w|^2} = \frac{(1+|\mu |)^2}{1-|\mu |^2}=K_{\mu}.
\end{equation}
Hence, it follows from \eqref{eq:4}, \eqref{eq:9} and our initial hypothesis $\cK |_{\cO} \in L_{\rm loc}^1 (\cO)$, see \eqref{co}-(iii), that $K(w,f)\in L^1$ in a neighborhood of the point $w_0=x_0+iy_0\in \C$. To see this, recall that by definition, $\cK=\Lambda/\landa$. Thus, we are in the conditions of the Iwaniec-Sverak theorem for degenerate elliptic quasiregular mappings (\cite{IS}, see also \cite{AIM}), which provides a Stoilow factorization for $f$ in a neighborhood of $w_0$. This implies that, around $w_0$, $f$ is either constant or an open mapping. We summarize this conclusion in the following assertion for later use:

\begin{assertion}\label{ass:1}
If $\nu_0 = \frac{1}{\sqrt{1+x_0^2+y_0^2}}(x_0,y_0,1)$ lies in $\cO\subset \S^2$, then $\nabla h$ is either an open mapping or constant around $(x_0,y_0)$.
\end{assertion}

{\bf Step 2:} \emph{Gradient mappings and support planes.}  

In Steps 2 through 9 of the proof of Theorem \ref{th:main}, we will let $u\in C^2(\R^3\setminus \{0\})$ be a degree one homogeneous solution to a linear equation \eqref{edp1}, and only assume that the coefficients $a_{ij}$ of \eqref{edp1} satisfy the degenerate ellipticity conditions (i), (ii) in \eqref{co}. That is, we will not use condition (iii) in \eqref{co}.

By homogeneity, $D^2u(x)$ always has a trivial zero eigenvalue corresponding to the radial direction, for any $x\in \R^3\setminus\{0\}$. Denote by $\mu_1(x)\leq \mu_2(x)$ the other two eigenvalues. These are also the eigenvalues of the spherical Hessian $\nabla_{\S^2}^2v$ of the function $v:=u(x/|x|)\in C^2(\S^2)$ at the point $\eta=x/|x|$, see e.g. \cite{GWZ}. Here, the spherical Hessian of $v$ is defined by $\nabla_{\S^2}^2v(\eta)=(v_{ij}(\eta) + v(\eta)\delta_{ij})$, where $v_{ij}$ are covariant derivatives with respect to a local orthonormal frame in $\S^2$. Then, the property that the coefficients $a_{ij}$ of \eqref{edp1} satisfy the degenerate ellipticity conditions i), ii) in \eqref{co} is equivalent to the fact that $\mu_1\mu_2\leq 0$ everywhere, i.e., to the fact that, on $\S^2$, ${\rm det}(\nabla_{\S^2}^2 v)\leq 0.$ This follows from the argument indicated after equation \eqref{eq:11}.

Consider the \emph{hedgehog} in $\R^3$ given by the restriction of the gradient mapping of $u$ to the unit sphere, $\nabla u : \S^2\flecha \R^3$. It can be regarded as a compact surface (with singularities) in $\R^3$, see \cite{LLR}. By compactness, $\nabla u(\S^2)$ admits a support plane in any direction, where here by a \emph{support plane} in the direction $\xi\in \S^2$ we mean a plane $\Pi_{\xi}\subset \R^3$ orthogonal to $\xi$ that touches $\nabla u(\S^2)$ at some point $q_{\xi}$, and so that $\esiz \nabla u-q_{\xi},\xi\esde\leq 0$ on $\S^2$.  Observe that $\nabla u(\S^2)$ cannot be constant, since $D^2 u$ is not identically zero. Thus, for almost every direction $\xi\in \S^2$, the two associated support planes to $\xi$ and $-\xi$ are different, and each of them intersects $\nabla u (\S^2)$ at a unique point.

Given arbitrary Euclidean coordinates $(x,y,z)$ in $\R^3$, the hedgehog $\nabla u:\S^2\flecha \R^3$ can be parametrized as the map in \eqref{eq:1g}, for all $\nu \in \S^2$ with positive $z$-coordinate, that is, 
 \begin{equation}\label{eq:9b}
\psi(x,y):=\nabla u (\nu)=(h_x,h_y,h-xh_x-yh_y), 
 \end{equation} 
where
\begin{equation}\label{eq:9bb}
 \nu:=\frac{(x,y,1)}{\sqrt{1+x^2+y^2}}.
\end{equation}
Recall that, by \eqref{eq:11}, $h_{xx} h_{yy}-h_{xy}^2\leq 0$. Obviously, $\psi(x,y)$ is an immersion with unit normal $\nu$ at the points where ${\rm det}(D^2 h)<0$. We call these points \emph{regular points} of the hedgehog. We should note that, although $\psi$ is at first only of class $C^1$, it can be easily checked using the inverse function theorem that any regular point $q$ of $\psi$ has a neighborhood $\cU\subset \R^2$ such that $\psi(\cU)$ is a $C^2$ graph over an open set of its tangent plane at $q$. Thus, it makes sense to talk about the second fundamental form $II$ of \eqref{eq:9b} at regular points, and a computation from \eqref{eq:9b}, \eqref{eq:9bb} shows that 
\begin{equation}\label{eq:10}
II=\frac{-1}{\sqrt{1+x^2+y^2}} D^2 h(x,y).
\end{equation}
In particular, the hedgehog has negative curvature at its regular points, and therefore such points cannot arise as contact points of $\nabla u(\S^2)$ with a support plane. Note that the hedgehog $\nabla u(\S^2)$ is regular at a point $\nu\in \S^2$ if and only if the two non-trivial eigenvalues $\mu_1\leq \mu_2$ of $D^2 u (\nu)$ are non-zero (and so, necessarily, of opposite signs), i.e. if and only if $D^2 u(\nu)$ has rank $2$.

\begin{definition}\label{pogo}
We say that $p_0\in \nabla u (\S^2)$ is a \emph{Pogorelov point} if there exists a direction $\xi\in \S^2$ such that $\nabla u (\S^2)\cap \Pi_{\xi}=\{p_0\}$, and $p_0 \not\in \{\nabla u(\xi), \nabla u (-\xi)\}$.
\end{definition}

\begin{assertion}\label{ass:2}
There exists a Pogorelov point of $\nabla u(\S^2)$. \end{assertion}
\begin{proof}
We first note that $\nabla u:\S^2\flecha \R^3$ has a regular point. Indeed, otherwise we would have $\mu_1\mu_2=0$ on $\S^2$. Thus, the function $f:=u|_{\S^2}$ would satisfy ${\rm det}(\nabla_{\S^2}^2 f )=0$ everywhere on $\S^2$. By \cite[Theorem 1]{K}, $f$ would be linear on $\S^2$. So, $u$ would also be linear, a contradiction.

Let then $\xi\in \S^2$ be a regular point of $\nabla u$. By slightly varying $\xi$, we can assume additionally that each of the support planes $\Pi_{\xi}$ and $\Pi_{-\xi}$ intersects $\nabla u(\S^2)$ at a unique point, say $q_1$ and $q_2$. As $\nabla u(\xi)$ cannot lie in any of these two planes (by regularity), either $q_1$ or $q_2$ is a Pogorelov point for $\nabla u(\S^2)$. 
\end{proof}

{\bf Step 3:} \emph{Setup for the rest of the proof.}

We fix from now on a Pogorelov point $p_0\in \nabla u(\S^2)$, with associated direction $\xi\in \S^2$. Take $\nu_0\in\S^2$ with $\nabla u(\nu_0)=p_0$. We consider Euclidean coordinates $(x,y,z)$ with $\xi=(1,0,0)$ and $\nu_0=(\nu_0^1,0,\nu_0^3)$, with $\nu_0^3>0$. One should observe that $\nu_0$ is not uniquely determined by $\xi$, since the subset $(\nabla u)^{-1}(p_0)$ of $\S^2$ might be large. As a matter of fact, we seek to show that it contains a geodesic semicircle. At this stage of the proof we will not require any additional information on $\nu_0$, but in Step 8 we will discuss how to choose it in a convenient way.

Since $\xi=(1,0,0)$, the support plane $\Pi_{\xi}$ leaves $\nabla u (\S^2)$ on its left side, i.e., $\Pi_{\xi}$ is of the form $x=\mu_{\rm max}$, and 
\begin{equation}\label{ja3}
\mu_{\rm min} \leq u_x(p) \leq \mu_{\rm max} \hspace{1cm} \forall p\in \S^2,
\end{equation}
for some values $\mu_{\rm min},\mu_{\rm max}\in \R$. The points $\nabla u (\pm \xi)$ do not lie in $x=\mu_{\rm max}$, since $p_0$ is a Pogorelov point. Thus, there exist $\mu_0<\mu_{\rm max}$ and $\ep>0$ such that $u_x(p)\leq \mu_0$ for every $p\in B(\xi;\ep)\cup B(-\xi; \ep)$, where here $B(a;\ep)$ denotes a geodesic ball in $\S^2$ of center $a$ and radius $\ep$. By homogeneity, $u_x(x,y,z)\leq \mu_0$ on a subset of $\R^3$ of the form $x^2\geq \delta (y^2+z^2)$ for some $\delta=\delta(\ep)>0$.

From now on, let $\Sigma$ be the entire saddle graph in $\R^3$ given by $z=h(x,y)$, where $h$ is defined by \eqref{eq:0}; note that $\Sigma$ has non-positive curvature at every point, by \eqref{eq:11}. By \eqref{eq:1g} and the compactness of $\nabla u (\S^2)$, we see that $\nabla h$ is uniformly bounded in $\R^2$. Moreover, by \eqref{ja3}, \eqref{eq:1g} and the definition of $\mu_0$, we have
 \begin{equation}\label{ja3b}
 \mu_{\rm min} \leq h_x(x,y) \leq \mu_{\rm max},
 \end{equation}
for all $(x,y)\in \R^2$, and 
 \begin{equation}\label{ja4}
 h_x(x,y) \leq \mu_0 <\mu_{\rm max}  \hspace{1cm} \text{$\forall (x,y)\in \R^2$ with $x^2 \geq \delta (y^2+1).$}
 \end{equation}
We will denote by $\Omega^+$ (for $x>0$) and $\Omega^-$ (for $x<0$) the two connected components of the set $x^2 \geq \delta (y^2+1)$ in $\R^2$. 
Also, note that 
\begin{equation}\label{ja4b}
h_x (x_0,0)=\mu_{\rm max}, \hspace{1cm} \text{ where } \nu_0=(\nu_0^1,0,\nu_0^3)=\frac{(x_0,0,1)}{\sqrt{1+x_0^2}}.
\end{equation}
We will use frequently in what follows the notation 
\begin{equation}\label{eq:fi}
\varphi(x,y):=(x,y,h(x,y)).
\end{equation}

{\bf Step 4:} \emph{A transverse line $L_n^*$ to $\Sigma\cap \{y=0\}$ with almost maximum slope.}

Consider a plane $\Pi$ given by $z=P(x,y):=ax+by+c$, with $a>\mu_0$. Then, for any $y_0\in \R$, we have by \eqref{ja4} and $a>\mu_0$ that the line $L_{y_0}\equiv \Pi\cap \{y=y_0\}$ is above (resp. below) the graph $z=h(x,y_0)$ as $x\to \8$ (resp. $x\to -\8$). In this way, there exist points $x_1(y_0)\leq x_2(y_0)$ such that
\begin{equation}\label{ja5}
h(x,y_0)>P(x,y_0)  \text{ for $x<x_1(y_0)$,} \hspace{1cm} h(x,y_0)< P(x,y_0)  \text{ for $x>x_2(y_0)$.}
\end{equation}
In particular, there exist points $(x_1,0)\in \Omega^-$ and $(x_2,0)\in \Omega^+$ such that $h(x,0)> P(x,0)$ for all $x\leq x_1$, and $h(x,0)< P(x,0)$ for all $x\geq x_2$. 

\begin{assertion}\label{ass:3}
There exist continuous curves $x=\alfa^-(y)$, $x=\alfa^+(y)$ in $\R^2$, which depend on the initial plane $\Pi$, such that $\alfa^-(0)=x_1$, $\alfa^+(0)=x_2$, and 
\begin{equation}\label{ja6}
h(\alfa^-(y),y)>P(\alfa^-(y),y), \hspace{1cm} h(\alfa^+(y),y)<P(\alfa^+(y),y) ,
\end{equation}
for all $y\in \R$.
\end{assertion}
\begin{proof}
Take $\bar{a}\in (\mu_0,a)$ and denote by $\overline{\mu}_{\rm min}$ the minimum value of $h_y$ in $\R^2$. Choose $\landa<0$ so that the half-line $\cL_{\landa}\subset \R^2$ given by $x=x_1+\landa y$ for $y\geq 0$ is contained in $\Omega^-$. We can obviously choose $\landa$ so that, additionally, $(a-\bar{a})\landa <\overline{\mu}_{\rm min} -b$ holds. See Figure \ref{fig:1}. Then, $\alfa^-(y):= x_1+\landa y$ satisfies the first inequality in \eqref{ja6} for all $y\geq 0$; indeed, if $(x,y)\in \cL_{\landa}$, integrating $\nabla h$ along $\cL_{\landa}$, and using that $h(x_1,0)>P(x_1,0)$ together with the previous inequalities we have 
$$h(x,y)
  > h(x_1,0) + (\bar{a}\landa + \overline{\mu}_{\rm min}) y \
  >  P(x_1,0)+ (a\landa +b)y =P(x,y).$$
The first inequality for $y<0$, and the second inequality in \eqref{ja6} are obtained similarly. This proves Assertion \ref{ass:3}.
\begin{figure}[htbp]
    \includegraphics[width=7.5cm]{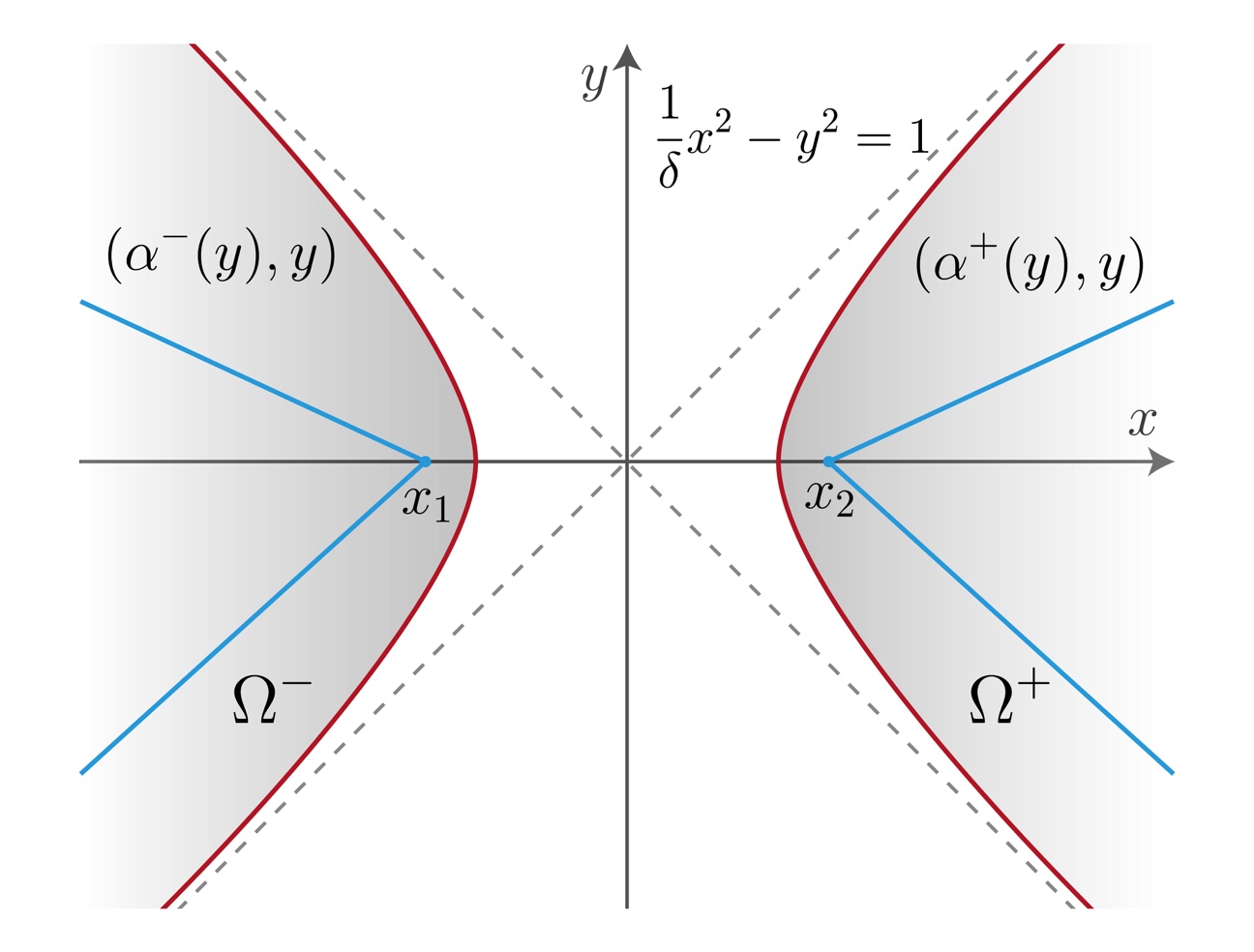}
     \caption{The curves $(\alfa^{\pm}(y),y)$ in $\R^2$.} 
\label{fig:1}
\end{figure}
\end{proof}

\begin{remark}\label{rem:alfa}
Observe that, if we consider the continuous curves $x=\alfa^{\pm} (y)$ defined in Assertion \ref{ass:3} with respect to the plane $\Pi$, then all points $\varphi(x,y)\in \Sigma$ where $x<\alfa^{-}(y)$ (resp. $x>\alfa^+(y)$) lie above (resp. below) $\Pi$. In order to see this, it suffices to realize that the proof of Assertion \ref{ass:3} also holds if, instead of $(x_1,0)\in \Omega^-$ we consider as initial point of $x=\alfa^-(y)$ any point $(x,0)$ with $x<x_1$ (and a similar argument for $x>x_2$ with $(x_2,0)\in \Omega^+$).
\end{remark}

Take next a sequence $\{\mu_n\}_n\to \mu_{\rm max}$, with $\mu_n\in (\mu_0, \mu_{\rm max})$ for all $n$. Consider the line $L_n$ in the vertical plane $y=0$ given by $z= \mu_n (x-x_0) + h(x_0,0)$. Note that $L_n$ intersects transversally $\Sigma_0:=\Sigma\cap \{y=0\}$ at $\varphi(x_0,0)$, by \eqref{ja4b}. More specifically, since $\mu_n< \mu_{\rm max}$, we see that $\Sigma_0$ lies below $L_n$ in the plane $y=0$ for values of $x<x_0$ near $x_0$, and above $L_n$ for $x>x_0$ near $x_0$. Besides, it is clear from \eqref{ja5} that $\Sigma_0$ lies above (resp. below) $L_n$ as $x\to -\8$ (resp. as $x\to \8$). This shows, in particular, that the planar set $\Sigma_0 \setminus L_n$ has at least four connected components, each of them homeomorphic to an open interval.

By the transversality of $\Sigma_0$ and $L_n$ at $\varphi(x_0,0)$, there exists some $\ep>0$ such that $h_x(x,0)>\mu_n$ and $\varphi(x,0)\not\in L_n$ for all $x\neq x_0$ with $|x-x_0|<\ep$. By Sard's theorem, if necessary, we can make a small parallel translation of $L_n$ in the plane $y=0$, to obtain a new straight line $L_n^*$ which might not pass through $(x_0,0,h(x_0,0))$ anymore, but which intersects $\Sigma_0$ transversely at every intersection point. Specifically, we may take $L_n^*$ so that it contains a point $\varphi(x_0^*,0)$ with $|x_0-x_0^*|<\ep$, and so that the distance between $\varphi(x_0^*,0)$ and $\varphi(x_0,0)$ is smaller than $1/n$. Here, $x_0^*=x_0^*(n)$, i.e., $x_0^*$ depends on $n$.

\begin{figure}[htbp]
    \includegraphics[width=7.8cm]{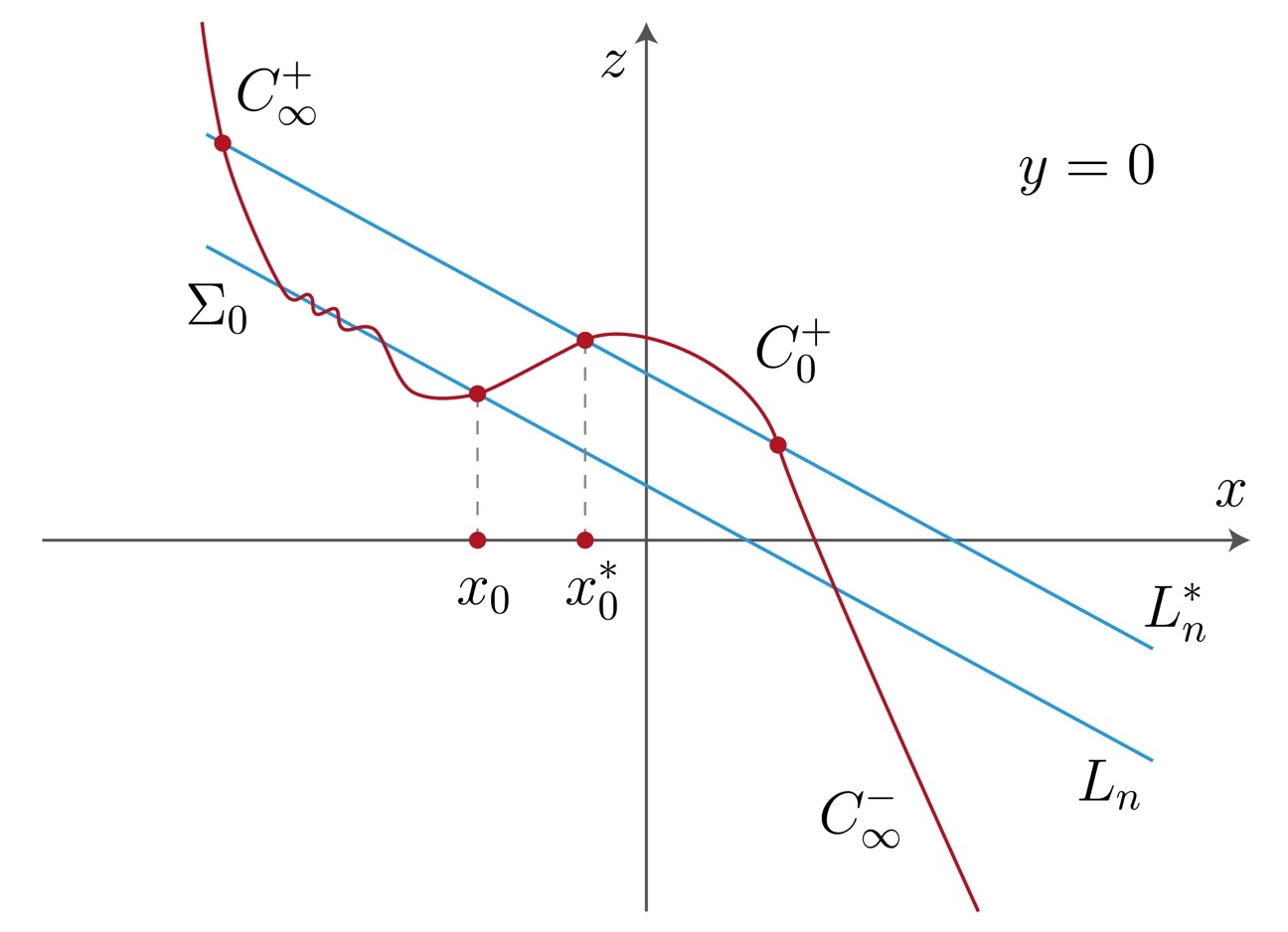}
     \caption{The connected components $C_{\8}^+$, $C_{\8}^-$ and $C_0^+$.} 
\label{fig:2}
\end{figure}

Note that, by \eqref{ja5}, $L_n^*$ lies either above or below $\Sigma_0$ as $x\to \8$ or $x\to -\8$. Then, by transversality, $\Sigma_0 \setminus L_n^*$ has a finite number of connected components. By the above arguments, we also know that the number of such connected components is at least $4$, and that $\varphi(x_0^*,0)$ lies at the common boundary of two such \emph{bounded} connected components. We will use the following notations for some special connected components of $\Sigma_0\setminus L_n^*$; (see Figure \ref{fig:2}).
\begin{enumerate}
\item
$C_{\8}^+$ is the unbounded component that lies strictly above $L_n^*$. 
 \item
$C_{\8}^-$ is the unbounded component that lies strictly below $L_n^*$.
 \item
$C_0^+$ is the \emph{bounded} component that lies strictly above $L_n^*$, and has $\varphi(x_0^*,0)$ as a boundary point.
\end{enumerate}
Observe that $C_{\8}^+$ lies in the set $\{x<x_0^*\}$, while $C_{\8}^-$ and $C_0^+$ lie in $\{x>x_0^*\}$

{\bf Step 5:} \emph{Study of the intersection of $\Sigma$ with the sheaf of planes containing $L_n^*$.}

Let us now fix the straight line $L_n^*$, and consider all the planes in $\R^3$, excluding $y=0$, that contain $L_n^*$. They are given by 
 \begin{equation}\label{plan}
z=P_b(x,y)= \mu_n (x-x_0^*) + b y + h(x_0^*,0),
 \end{equation} 
for each $b\in \R$. Call $\Pi_b$ to the plane determined by $b$. We next study $\Sigma\cap \Pi_b$.

Fix some point $q_0\in C_0^+$. Let $I_n$ (resp. $J_n$) denote the set of values $b\in \R$ for which $q_0$ can be joined to a point $\varphi(x,y)\in \Sigma$, with $y>n$ (resp. with $y<-n$), through an arc contained in $\Sigma\setminus (\Pi_b \cup C_{\8}^+).$ The statement of the next assertion uses that  $\overline{\mu}_{\rm min} \leq h_y (x,y)\leq \overline{\mu}_{\rm max}$ for adequate constants, for all $(x,y)\in \R^2$. It states that for any $n\in \N$ there exists a plane $\Pi_{b_n}$ such that we can find an arc in $\Sigma$ joining $q_0$ to points $\varphi(x,y)$ where $y>n$ and $y<-n$, while avoiding both $\Pi_{b_n}$ and the connected component $C_{\8}^+$.

\begin{assertion}\label{lem2}
There exists $b_n\in I_n\cap J_n$, with $\overline{\mu}_{\rm min}\leq b_n\leq \overline{\mu}_{\rm max}$.
\end{assertion}
\begin{proof}
Write $q_0=\varphi(q_0^1,0)$. By construction, $q_0$ lies above $L_n^*$. If we choose $b\leq \overline{\mu}_{\rm min}$, then $\varphi(q_0^1,y)\in \Sigma$ lies above $\Pi_b$, for all $y>0$. Since $\varphi(q_0^1,0)\not\in C_{\8}^+$, this means that $\overline{\mu}_{\rm min}\in I_n$. By the same argument, $\overline{\mu}_{\rm max}\in J_n$. Thus $I_n$ and $J_n$ are non-empty, and they both intersect the closed interval $[\overline{\mu}_{\rm min},\overline{\mu}_{\rm max}]$.

We check next that $I_n$ is open. Let $b_0\in I_n$. Then, there exists an arc in $\Sigma\setminus (\Pi_{b_0}\cup C_{\8}^+)$ joining $q_0$ with a point $p=\varphi(x,y)$, with $y>n$. By compactness, this arc lies above $\Pi_{b_0}$ at a certain distance $d>0$. In particular, for values of $b$ near $b_0$, this arc also avoids $\Pi_{b}\cup C_{\8}^+$. Therefore, $I_n$ is open. By the same argument, $J_n$ is open.

Finally, we prove that $I_n\cup J_n=\R$, what, together with the already proved properties and the fact that $[\overline{\mu}_{\rm min},\overline{\mu}_{\rm max}]$ is connected, yields Assertion \ref{lem2}. Arguing by contradiction, assume that there exists $b\in \R\setminus (I_n\cup J_n)$. We are going to prove next that the (open) connected component of $\Sigma\setminus \Pi_b$ that contains $q_0$, which we will denote by $\Sigma (C_0^+)$, is \emph{bounded}. This will contradict the fact that $\Sigma$ is a saddle graph.

To do this, we start fixing some notation and making some elementary comments. First, note that $\Sigma (C_0^+)$ lies above $\Pi_b$, since $q_0\in C_0^+$. Also, denote by $\Sigma(C_{\8}^+)$ the connected component of $\Sigma\setminus \Pi_b$ that contains $C_{\8}^+$. By Remark \ref{rem:alfa}, if we consider the continuous curves $x=\alfa^{\pm} (y)$ defined in Assertion \ref{ass:3} with respect to the plane $\Pi_b$, then all points $\varphi(x,y)\in \Sigma$ where $x<\alfa^{-}(y)$ (resp. $x>\alfa^+(y)$) lie above (resp. below) $\Pi_b$. In this way, the curve $\Gamma^-:=\{\varphi(\alfa^{-}(y),y)\in\Sigma : y\in \R\}$ is contained in $\Sigma(C_{\8}^+)$.

First of all, we prove that every point $\varphi(x,y)$ of $\Sigma(C_0^+)$ satisfies $y\in [-n,n]$. Indeed, otherwise, there would exist an arc $\gamma$ in $\Sigma$ starting at $q_0$, that reaches either $\{y<-n\}$ or $\{y>n\}$, and that intersects $C_{\8}^+$, since $b\not\in I_n\cap J_n$. Let $\overline{z}_0$ denote the first point where $\gamma$ touches $C_{\8}^+$. Then, a neighborhood of $\overline{z}_0$ trivially lies in $\Sigma(C_{\8}^+)$. See Figure \ref{fig:compo}. In particular, $\Sigma(C_{\8}^+)= \Sigma(C_0^+)$. 

\begin{figure}[htbp]
     \includegraphics[width=10cm]{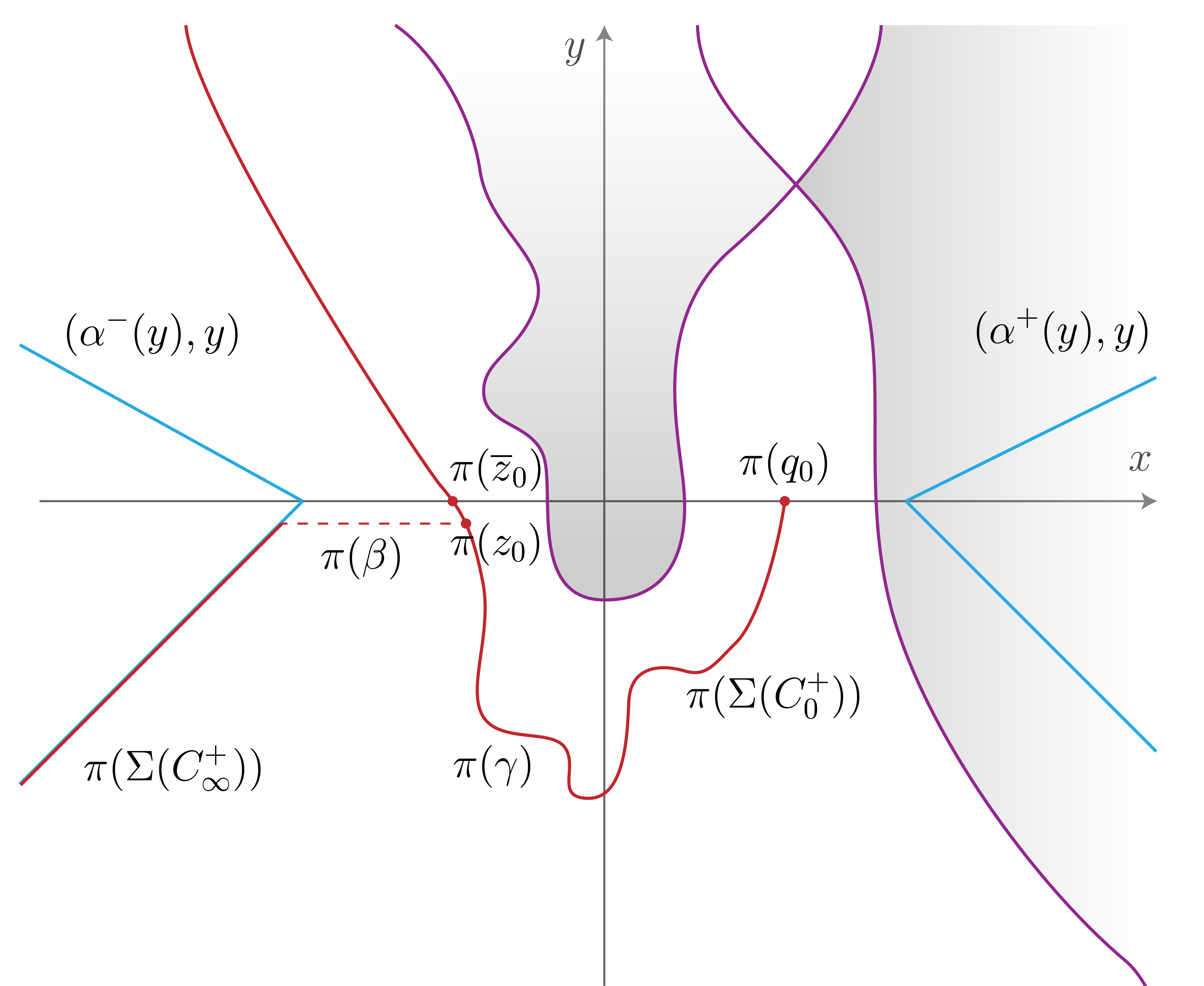}
     \caption{Proof that $\Sigma(C_0^+)$ lies in the slab of $\R^3$ given by $|y|\leq n$. In the figure, $\pi$ denotes the projection $\pi(x,y,z)=(x,y)$ onto the $x,y$-plane}.
\label{fig:compo}
\end{figure}

Let $z_0=\varphi(z_0^1,z_0^2)$ be a point of that neighborhood, that also lies in the interior of the arc of $\gamma$ between $q_0$ and $\overline{z}_0$. Assume that $z_0^2<0$ (the argument is similar if $z_0^2>0$). Then, we can join the curve $\Gamma^-\subset\Sigma$ defined above with the point $z_0$ along an arc $\beta$ contained in $\Sigma(C_{\8}^+)$ and so that every point of the arc $\beta$ has negative $y$-coordinate. See Figure \ref{fig:compo}. This implies that $\beta$ does not touch $C_{\8}^+$, which is contained in the $y=0$ plane. Now, the union of the arc of $\gamma$ joining $q_0$ with $z_0$, the arc $\beta$, and a suitable arc of the curve $\Gamma^-$ produces an arc in $\Sigma(C_{\8}^+)$ that avoids $C_{\8}^+$ and joins $q_0$ with a point in $\Sigma\cap \{y<-n\}$ (see Figure \ref{fig:compo}). This would mean that $b\in J_n$, a contradiction. Thus, $\Sigma(C_0^+)$ lies in the slab of $\R^3$ given by $|y|\leq n$, as desired.

Recall that all points of the form $\varphi(\alfa^+(y),y)$ lie below $\Pi_b$, by Assertion \ref{ass:3}. Since all points $\varphi(x,y)\in \Sigma(C_0^+)$ satisfy $|y|\leq n$ and lie above $\Pi_b$, we conclude then that their $x$-coordinates are bounded from above by $\alfa^+(y)$.

On the other hand, assume that there exists an arc in $\Sigma (C_0^+)$ that joins $q_0$ with a point of the form $\varphi(\alfa^-(y),y)$. By Assertion \ref{ass:3}, we have $\Sigma(C_0^+)=\Sigma(C_{\8}^+)$.  But now, as $\Sigma(C_{\8}^+)$ has points of the form $\varphi(x,y)$ with $|y|$ arbitrarily large, we contradict the fact that $\Sigma(C_0^+)$ lies in the slab $|y|\leq n$. 

We have then proved that $\Sigma(C_0^+)$ is contained in the compact set $$\{\varphi(x,y): \alfa^-(y)\leq x \leq \alfa^+(y), |y|\leq n\}\subset \Sigma.$$ Thus $\Sigma(C_0^+)$ is a bounded connected component of $\Sigma\setminus \Pi_b$, in contradiction with the saddleness of $\Sigma$. This proves Assertion \ref{lem2}.
\end{proof}

{\bf Step 6:} \emph{Study of the intersection of $\Sigma$ with the limit plane $\Pi_{\8}$.}

\vspace{0.2cm}

For each $n$, let $b_n \in \R$ be given by Assertion \ref{lem2}, and consider the associated plane $\Pi_{b_n}$ given by \eqref{plan} for $b=b_n$. Since $\overline{\mu}_{\rm min}\leq b_n\leq \overline{\mu}_{\rm max}$, we have up to subsequence that $\{b_n\}_n\to b_{\8}\in [\overline{\mu}_{\rm min},\overline{\mu}_{\rm max}]$. Since $|x_0^*-x_0|<1/n$ and $\{\mu_n\}_n\to \mu_{\rm max}$, the planes $\Pi_{b_n}$ converge to the limit plane 
\begin{equation}\label{limplan}
\Pi_{\8} \hspace{0.2cm} : \hspace{0.2cm} z=P_{\8}(x,y):= \mu_{\rm max} (x-x_0) + b_{\8} y + h(x_0,0),
\end{equation}
which passes through $\varphi(x_0,0)\in \Sigma$ with maximum slope $\mu_{\rm max}$ in the $x$-direction. 

We study next $\Sigma\cap \Pi_{\8}$. Fix any $y_0\in \R$. Then, taking $\Pi=\Pi_{\8}$ in Assertion \ref{ass:3}, it is a consequence of \eqref{ja6} that the curve $\Sigma\cap \{y=y_0\}$ intersects $\Pi_{\8}$. 

\begin{assertion}\label{lem3}
Either for all $y_0\geq 0$, or for all $y_0\leq 0$, there exist $x^-(y_0)\leq x^+(y_0)$ such that
$$\Pi_{\8}\cap \Sigma\cap \{y=y_0\}=\{\varphi(x,y_0): x\in [x^-(y_0),x^+(y_0)]\}.$$ Moreover, $h_x(x,y_0)=\mu_{\rm max}$ holds for every $x\in [x^-(y_0),x^+(y_0)]$, and $\Sigma$ lies above $\Pi_{\8}$ (resp. below $\Pi_{\8}$) when $x<x^-(y_0)$ (resp. $x>x^+(y_0)$).
\end{assertion}
\begin{proof}
Fix $y_0\in \R$. We distinguish two possible situations.

\vspace{0.2cm}

{\bf Case 1:}\emph{ $\Pi_{\8}\cap \Sigma\cap \{y=y_0\}$ is not a unique point.} In that case, given two points $\varphi(x_1,y_0)$, $\varphi(x_2,y_0)$ of that intersection, we have that all points of the form $\varphi(x,y_0)$ with $x\in [x_1,x_2]$ also lie in $\Pi_{\8}\cap \Sigma\cap \{y=y_0\}$. This follows since $h_x(x,y_0)\leq \mu_{\rm max} = (P_{\8})_x(x,y_0)$ and $h(x_i,y_0)=P_{\8}(x_i,y_0)$, for $i=1,2$.
Thus, if $\Pi_{\8}\cap \Sigma\cap \{y=y_0\}$ has at least two points, there exist $x^-(y_0)<x^+(y_0)$ such that:

\begin{enumerate}
\item
$\varphi(x,y_0)$ lies above $\Pi_{\8}$ for all $x<x^-(y_0)$.
 \item
$\varphi(x,y_0)$ lies below $\Pi_{\8}$ for all $x>x^+(y_0)$.
 \item
$\varphi(x,y_0)\in\Pi_{\8}$ for all $x\in [x^-(y_0),x^+(y_0)]$.
\end{enumerate}
Note that $h_x(x,y_0)=\mu_{\rm max}$ for all $(x,y_0)$ in the third situation above. So, the statement of Assertion \ref{lem3} holds for every $y_0\in \R$ such that $\Pi_{\8}\cap \Sigma\cap \{y=y_0\}$ is not a unique point. No sign assumption is needed here for $y_0$.

\vspace{0.2cm}

{\bf Case 2:} \emph{$\Pi_{\8}\cap \Sigma\cap \{y=y_0\}$ is a unique point.} This situation is subtler, and needs an additional control on the intersections $\Sigma\cap \Pi_{b_n}$ before passing to the limit. 

Let $b_n$ be given by Assertion \ref{lem2}, with $\{b_n\}_n\to b_{\8}$. 
By $b_n \in I_n\cap J_n$, there exists an arc $\gamma^+=\gamma^+(n)$ in $\Sigma$ that lies above $\Pi_{b_n}$, that does not intersect $C_{\8}^+$, and whose endpoints have $y$-coordinate equal to $n$ and $-n$, respectively. Since $L_n^*$ intersects $\Sigma_0=\Sigma\cap \{y=0\}$ transversely at a finite number of points, there obviously exists a unique connected component $C_1^-$ of $\Sigma_0\setminus L_n^*$ that has as a boundary point the unique boundary point of $C_{\8}^+$, and lies \emph{below} $\Pi_{b_n}$  (since $C_1^-$ lies below $L_n^*$). As $\gamma^+$ lies above $\Pi_{b_n}$ and does not intersect $C_{\8}^+$, we easily deduce that every point in $\gamma^+\cap \Sigma_0$ is of the form $\varphi(x,0)$, with $x>{\rm sup} \{ x: \varphi(x,0)\in C_1^-\}$. Obviously, $\gamma^+\cap \Sigma_0$ is non-empty since $\gamma^+$ goes from $y=n$ to $y=-n$.

Let $\Sigma(C_1^-)$ denote the connected component of $\Sigma\setminus \Pi_{b_n}$ that contains $C_1^-$ (thus, it lies below $\Pi_{b_n}$). For each $n$, let $\alfa_n^+(y)$, $\alfa_n^-(y)$ be the functions $\alfa^+(y), \alfa^-(y)$ defined by Assertion \ref{ass:3} with respect to $\Pi=\Pi_{b_n}$. 
Then, $\Sigma(C_1^-)$ must intersect either $\Sigma\cap \{y=n\}$ or $\Sigma\cap \{y=-n\}$; indeed, otherwise, $\Sigma(C_1^-)$ would be a connected component contained in a compact region of $\Sigma$ bounded by $\gamma^+$, $\Sigma\cap \{y=\pm n\}$ and $\{\varphi(\alfa_n^{-}(y),y): y\in \R\}$, and this contradicts the saddleness of $\Sigma$. 

In this way, we can take an arc $\gamma^-=\gamma^-(n)$ contained in $\Sigma(C_1^-)$ that joins a point of $C_1^-$ with a point $q_n$ of $\Sigma$ with $y$-coordinate equal to $n$ or $-n$. Up to a subsequence of the $\{b_n\}_n$, we can assume that one of these two situations holds for all $n$. For definiteness, \emph{we will assume that the $y$-coordinate of $q_n$ is equal to $n$, for all $n$.}

Then, obviously, any plane $\{y=y_0\}$ with $y_0\in [0,n]$ is intersected by the curves $\gamma^-$, $\gamma^+$, and $\{\varphi(\alfa_n^+(y),y): y\in \R\}$. Using again that $\gamma^+\cap C_{\8}^+ = \emptyset$, we deduce the existence of points $x_1<x_2<\alfa_n^+(y_0)$, with each $x_1,x_2$ depending on $y_0$ and $n$, such that $$\varphi(x_1,y_0)\in \gamma^-,\hspace{0.5cm} \varphi(x_2,y_0)\in \gamma^+.$$
Therefore, there exist $x_3\in (x_1,x_2)$ and $x_4\in (x_2,\alfa^+(y))$ such that both $\varphi(x_3,y_0)$ and $\varphi (x_4,y_0)$ lie in $\Sigma\cap\Pi_{b_n}\cap \{y=y_0\}$. Besides, since the line $\Pi_{b_n}\cap \{y=y_0\}$ has slope $\mu_n$ and $\varphi(x_2,y_0)$ lies above $\Pi_{b_n}$, with $x_2\in (x_3,x_4)$, by the mean value theorem there must exist $x_5\in (x_3,x_4)$ such that $\varphi(x_5,y_0)$ lies above $\Pi_{b_n}$, and $h_x(x_5,y_0)=\mu_n$.

From now on, we denote $s_n(y_0):=x_3<t_n(y_0):=x_5$. Thus, for every $n\in \N$ and every $y\in [0,n]$, we have:
\begin{enumerate}
\item
$\varphi(s_n(y),y)\in \Sigma\cap \Pi_{b_n}$.
 \item
$\varphi(t_n(y),y)$ lies above $\Pi_{b_n}$, and $h_x(t_n(y),y)=\mu_n$.
\end{enumerate}

We now pass to the limit, and show that the statement of Assertion \ref{lem3} holds for every $y_0\geq 0$; if we had assumed that the $y$-coordinate of $q_n$ is $-n$, the next argument would show that Assertion \ref{lem3} holds for every $y_0\leq 0$.

Fix then $y_0\geq 0$. By our hypothesis in the present Case 2 and \eqref{ja5}, there exists a certain value $x(y_0)$ such that $\varphi(x,y_0)$ lies above $\Pi_{\8}$ for all $x<x(y_0)$, and below $\Pi_{\8}$ for all $x>x(y_0)$. 

Take $(c(y_0),y_0)\in \Omega^-$ with $c(y_0)<x(y_0)$. Since $\{\Pi_{b_n}\}_n\to \Pi_{\8}$, there exists $n_0\in \N$ such that 
$\varphi(c(y_0),y_0)$ lies above $\Pi_{b_n}$, for every $n\geq n_0$. Now, as $(c(y_0),y_0)\in \Omega^-$, we have by \eqref{ja4} and $\mu_0<\mu_n$ that $\varphi(x,y_0)$ lies above $\Pi_{b_n}$, for all $x<c(y_0)$ and all $n\geq n_0$. In particular, $c(y_0)<s_n(y_0)<t_n(y_0)$, for all $n$ large enough, since $\varphi(s_n(y_0),y_0)\in \Pi_{b_n}$.

Arguing in a similar way for large positive values of $x$, we deduce that the sequences $\{s_n(y_0)\}_n$ and $\{t_n(y_0)\}_n$ are bounded. Thus, up to a subsequence, we must have $\{\varphi(s_n(y_0),y_0)\}_n\to \varphi(x(y_0),y_0)$, by uniqueness of the point $\varphi(x(y_0),y_0)$.

On the other hand, the points $\varphi(t_n(y_0),y_0)$ converge to some point that is not below $\Pi_{\8}$, since $\varphi(t_n(y_0),y_0)$ lies above $\Pi_{b_n}$ and $\{\Pi_{b_n}\}_n\to \Pi_{\8}$. But since $t_n(y_0)>s_n(y_0)\to x(y_0)$ and $\varphi(x,y_0)$ lies below $\Pi_{\8}$ for all $x>x(y_0)$, we deduce then that $\{t_n(y_0)\}_n\to x(y_0)$. In particular, $h_x(x(y_0),y_0) = \mu_{\rm max}$, since $h_x(t_n(y_0),y_0)=\mu_n$. This proves Assertion \ref{lem3} in Case 2, and thus completes the proof. 
\end{proof}

{\bf Step 7:} \emph{Existence of a half-line of maximal slope in $\Sigma\cap \Pi_{\8}$.} 

\vspace{0.2cm}

In this step, we show that the set $\Sigma\cap \Pi_{\8}$ contains some half-line $\cL^*$, and moreover, $h_x (x,y)= \mu_{\rm max}$ for all $(x,y)\in \R^2$ with $\varphi(x,y)\in \cL^*$.

To start, assume for definiteness that Assertion \ref{lem3} holds for $y_0\geq 0$ (the case $y_0\leq 0$ is treated analogously). Let $\mathcal{J}$ be the set of values $y_0\in \R$ such that $\Pi_{\8}\cap \Sigma\cap \{y=y_0\}$ is a unique point $\varphi(x(y_0),y_0)$, at which $h_x(x(y_0),y_0)<\mu_{\rm max}$ holds. Then, by Assertion \ref{lem3}, we have $\mathcal{J}\subset (-\8,0)$. Let $\delta_0\leq 0$ denote the supremum of $\mathcal{J}$, where we use the convention that $\delta_0=-\8$ if $\mathcal{J}$ is empty.

It follows from Assertion \ref{lem3} that there exist two (at first, maybe non-continuous) functions $x^-(y)<x^+(y)$, defined for all $y>\delta_0$, and such that the following properties hold:
\begin{equation}\label{i.ii}
\left\{\def\arraystretch{1.5}\begin{array}{l}
\text{ $i)$ \hspace{0.1cm} $h(x,y)>P_{\8}(x,y)$  if $x<x^-(y)$.}\\
\text{ $ii)$ \hspace{0.1cm}$h(x,y)<P_{\8}(x,y)$  if  $x>x^+(y)$.} \\
\text{ $iii)$ \hspace{0.1cm}$h(x,y)= P_{\8}(x,y)$ and  $h_x(x,y)=\mu_{\rm max}$  if $x\in [x^-(y),x^+(y)]$.}
\end{array}\right.
\end{equation}
To see this, one should recall that our conclusion in Case 1 in the proof of Assertion \ref{lem3} holds for all $y_0\in \R$, not only for $y_0\geq 0$ or $y_0\leq 0$.

\begin{assertion}\label{lem4}
The sets $$D^-=\{(x,y)\in \R\times (\delta_0,\8) : h(x,y)>P_{\8}(x,y)\},$$ $$D^+=\{(x,y)\in \R\times (\delta_0,\8) : h(x,y)<P_{\8}(x,y)\}$$ are open convex sets of $\R^2$. In particular, $x^+(y)$, $x^-(y)$ are continuous.
\end{assertion}
\begin{proof}
We will prove the result for $D^+$; the argument for $D^-$ is analogous. Let $p_i:=(x_i,y_i)\in D^+$, $i=1,2$. If $y_1=y_2$, the segment that joins both points lies in $D^+$, by property $ii)$ in \eqref{i.ii}.

Assume that $y_1\neq y_2$, and that the segment that joins $p_1$ with $p_2$ is not contained in $D^+$. As $x^+(y)<\alfa^+(y)$ and $\alfa^+(y)$ is continuous, we can take a translation of $\overline{p_1p_2}$ in the positive $x$-direction so that the resulting segment is contained in $D^+$. Next, translate that segment back in the negative $x$-direction, until reaching a first contact point with the set $D_0:=\{(x,y): h(x,y)=P_{\8}(x,y)\}$. We will denote the resulting segment by $\cS_0$. 

Note that the endpoints of $\cS_0$ lie in $D^+$, and that $D^+$ is connected by properties $i)$-$iii)$ in \eqref{i.ii}. Let $\gamma$ denote a compact arc in $D^+$ joining the endpoints of $\cS_0$. Then, there exists $\varepsilon>0$ such that $h\leq P_{\8}-\varepsilon$ for any point of $\gamma$. In this way, if we let $r_{\8}$ denote the line in the intersection of $\Pi_{\8}$ with the vertical plane that projects over the segment $\cS_0$, since $h\leq P_{\8}$ along $\cS_0$, we obtain the existence of a plane $\Pi_1$ that contains $r_{\8}$, has slope smaller than $\mu_{\rm max}$ in the $x$-direction, and does not touch $\varphi(\gamma)$; see Figure \ref{fig:convex}.

\begin{figure}[htbp]
    \includegraphics[width=11cm]{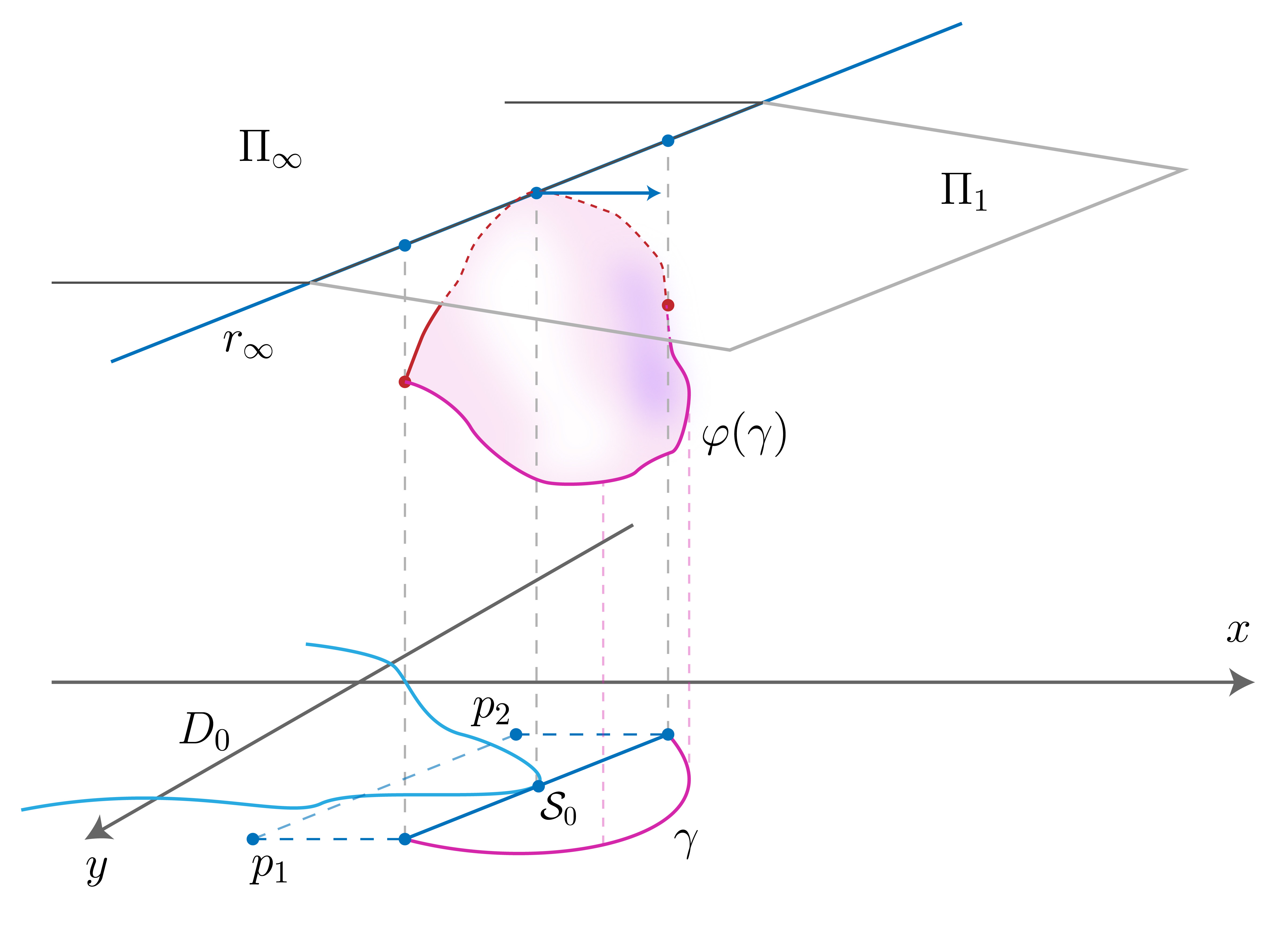}
     \caption{The argument in the proof of Assertion \ref{lem4}.} 
\label{fig:convex}
\end{figure}

Consider next the graph $\cG$ in $\R^3$ given by the restriction of $z=h(x,y)$ to the compact domain of $\R^2$ bounded by the segment $\cS_0$ and the curve $\gamma$. Since $\cG$ is saddle and its boundary does not touch the half-space of $\R^3$ above $\Pi_1$, then $\cG$ also has this property. But now, observe that at the points of the non-empty set $\cS_0\cap D_0$ we have $h_x=\mu_{\rm max}$. Since the slope of $\Pi_1$ in the $x$-direction is smaller than $\mu_{\rm max}$, this implies that there should exist points of $\cG$ above $\Pi_1$, a contradiction. This proves Assertion \ref{lem4}.
\end{proof}

Since $D^-$, $D^+$ are disjoint, open convex sets of $\R^2$, there exists a line $\cL\subset \R^2$ that separates them strictly, i.e., $D^-$ and $D^+$ lie in different connected components of $\R^2-\cL$. In particular, any point of the straight half-line $\cL^*:=\cL\cap \{y\geq \delta_0\}$ lies in the set 
 \begin{equation}\label{eq:dee}
\cD=\{(x,y): y \geq \delta_0, x\in [x^-(y),x^+(y)]\}.
\end{equation} 
Observe that, by iii) of \eqref{i.ii}, we have $h_x=\mu_{\rm max}$ and $h=P_{\8}$ on $\cD$, i.e., $\varphi(\cD)\subset \Pi_{\8}\cap \Sigma$. Since the intersection of $\nabla u(\S^2)$ with the support plane $x=\mu_{\rm max}$ of $\R^3$ is just the point $p_0$, we deduce that $\psi(\cD)=\{p_0\}$, where $\psi$ is given by \eqref{eq:9b}. Thus, $h_y$ is constant on $\cD$. In particular, $h_x$ and $h_y$ are constant along $\cL^*$, with $h_x=\mu_{\rm max}$. Then, $\varphi (\cL^*)$ is a straight half-line that lies in $\Sigma\cap \Pi_{\8}$, and we deduce from there that $h_y=b_{\8}$ on $\cD$, where $b_{\8}$ is defined in \eqref{limplan}. In particular, the limit plane $\Pi_{\8}$ is tangent to $\Sigma$ at every point of $\varphi(\cD)$. Also,
\begin{equation}\label{pecero}
p_0= (\mu_{\rm max}, b_{\8}, *)\in \R^3.
\end{equation}
Note that if $\delta_0=-\8$, both $\cL^*$ and $\varphi(\cL^*)$ are (complete) lines.

{\bf Step 8:} \emph{Existence of a geodesic semicircle in $(\nabla u)^{-1}(p_0)$.} 

In this step we show that, by choosing in a more careful way the initial direction $\nu_0\in (\nabla u)^{-1}(p_0)$ that we fixed at the beginning of Step 3, we can ensure that $\Omega_{\xi}:=(\nabla u)^{-1}(p_0)$ contains a geodesic semicircle of $\S^2$. 

\begin{figure}[htbp]
    \includegraphics[width=7cm]{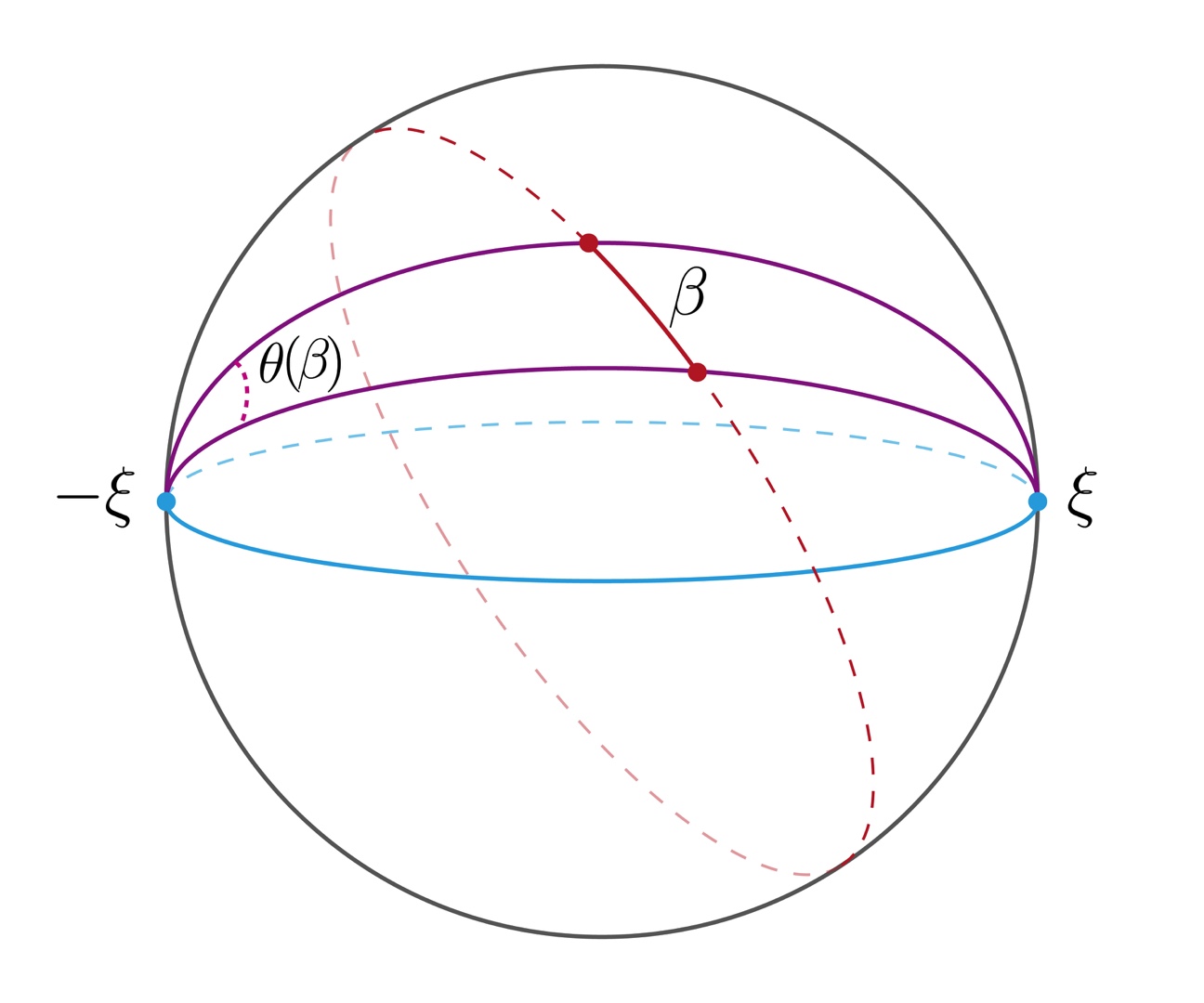}
     \caption{The definition of angle $\theta(\beta)$.} 
\label{fig:angledef}
\end{figure}

Assume that this last property is not true. Let $\beta$ be any geodesic arc of $\S^2$ contained in $\Omega_{\xi}$, and denote its endpoints by $\{\beta_0^1,\beta_0^2\}$. Note that, by our choice of the direction $\xi$ in Step 3, the distance in $\S^2$ between the compact subsets $\Omega_{\xi}$ and $\{\xi,-\xi\}$ is positive (since $p_0$ is a Pogorelov point). Thus, we can consider the \emph{angle} $\theta(\beta)\in [0,\pi]$ at $\xi$ defined by the two geodesic semicircles $\gamma_1,\gamma_2$ of $\S^2$ with endpoints $\{\xi,-\xi\}$ that satisfy $\beta_0^i \in \gamma_i$. See Figure \ref{fig:angledef}. Since $\beta$ has length $<\pi$ by hypothesis, this angle is $<\pi$. 

Observe first of all that there exists at least one geodesic arc (of positive length) $\beta^*$ contained in $\Omega_{\xi}$. To see this, let $\cL^*$ denote the straight half-line of the $(x,y)$-plane whose existence was shown in Step 7. Let $\beta^*$ be the geodesic arc in $\S^2$ that corresponds to $\cL^*$ via the totally geodesic bijection $\R^2\flecha \S^2_+$ given by \eqref{eq:9bb}. Since $h_x=\mu_{\rm max}$ along $\cL^*$, we have from \eqref{eq:9b} and \eqref{pecero} that 
\begin{equation}\label{eq:betastar}
\beta^*\subset (\nabla u)^{-1}(p_0)=\Omega_{\xi}.
\end{equation}
Since $\cL^*$ is not parallel to the $y$-axis, clearly $\theta(\beta^*)>0$.

We next prove that there exists a geodesic arc $\beta_{\8}$ of maximum angle in $\Omega_{\xi}$. Let $\theta_0\in (0,\pi]$ denote the supremum of the angles $\theta(\beta)$, among all possible choices of geodesic arcs $\beta$ contained in $\Omega_{\xi}$. Take any sequence of geodesic arcs $\{\beta_n\}_n$ in $\Omega_{\xi}$ with $\theta(\beta_n)\to \theta_0$. Then, up to a subsequence, the endpoints $a_n,b_n$ and the midpoint $c_n$ of the $\beta_n$ converge to three geodesically aligned points $\{a_1,a_2,a_3\}$ in $\Omega_{\xi}$. Since any point of $\beta_n$ is a convex combination of its endpoints, we deduce that $\{\beta_n\}_n$ converges  to the geodesic arc $\beta_{\8}$ contained in $\Omega_{\xi}$ with endpoints $\{a_1,a_2\}$ and midpoint $a_3$. In particular, $\beta_{\8}$ has positive length $<\pi$, and $\theta(\beta_{\8})=\theta_0$. We then conclude that $\theta_0<\pi$. 

\begin{figure}[htbp]
    \includegraphics[width=9.5cm]{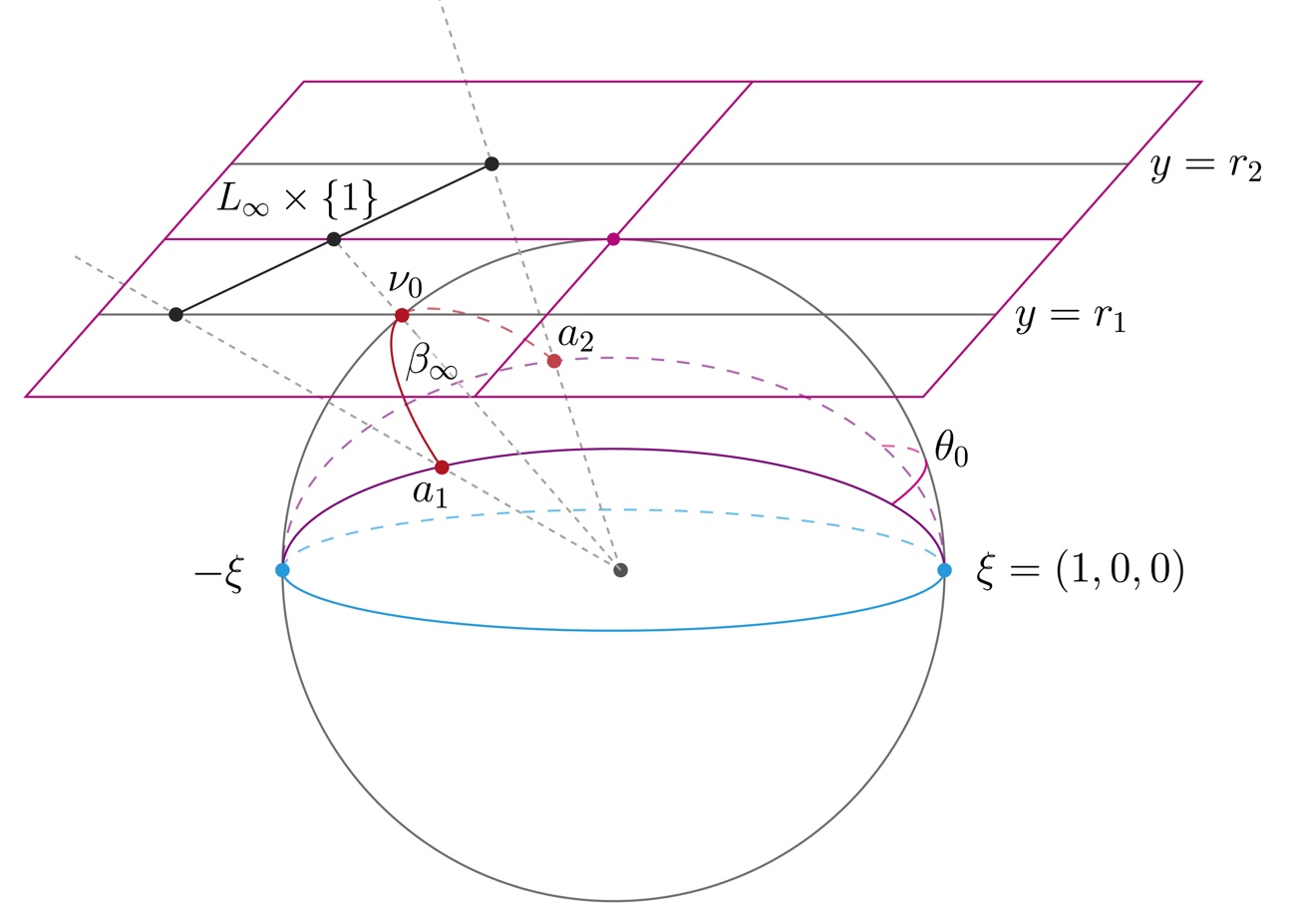}
     \caption{Choice of $\nu_0\in (\nabla u )^{-1}(p_0)$.} 
\label{fig:angle}
\end{figure}

Once we know this property, it is clear that we can choose the original $\nu_0\in (\nabla u )^{-1}(p_0)$, which was initially chosen in Step 3 without any a priori limitation, as follows:  \emph{$\nu_0$ is the unique point of the geodesic arc $\beta_{\8}\subset \Omega_{\xi}$ with the property that the angles $\theta_1,\theta_2$ of the two geodesic arcs of $\beta_{\8}$ joining $\nu_0$ with each of the endpoints $\{a_1,a_2\}$ of $\beta_{\8}$ satisfy $\theta_i=\theta_0/2<\pi/2$, for $i=1,2$.} See Figure \ref{fig:angle}.
This choice for $\nu_0$ lets us choose in a more specific way the coordinates $(x,y,z)$ at the beginning of Step 3. Recall that, in these $(x,y,z)$ coordinates, we had $\xi=(1,0,0)$,  $\nu_0=(\nu_0^1,0,\nu_0^3)$ with $\nu_0^3>0$. By our new specific choice of $\nu_0$, after a suitable rotation of the $(x,y,z)$-coordinates around the $x$-axis, we can additionally suppose that the arc $\beta_{\8}$ lies in the hemisphere $\S^2\cap \{z>0\}$. Note that $\nu_0\in \beta_{\8}$, and that every point of $\beta_{\8}$ lies in $(\nabla u)^{-1}(p_0)$. 

Consider the totally geodesic bijection $\R^2\flecha \S^2_+$ given by \eqref{eq:9bb}. This bijection takes $\nu_0$ to $(x_0,0)$ for some $x_0\in \R$, and $\beta_{\8}$ to a compact line segment $L_{\8}$ passing through $(x_0,0)$. See Figure \ref{fig:angle}. In the same way, the geodesic semicircles $\gamma_1,\gamma_2$ in $\S^2\cap \{z\geq 0\}$ that pass through the points $\{\xi,-\xi,a_i\}$ are projected into two parallel lines in $\R^2$ of the form $y=r_i$, for some $r_1<0<r_2$. Obviously, each of the endpoints of $L_{\8}$ lies in one of these lines. 

We can now carry out the argument in Steps 3 through 7 for this new choice of $\nu_0$. Let $\cD\subset \R^2\cap \{y\geq \delta_0\}$ denote the subset given by \eqref{eq:dee} in Step 7 of the proof. Since $\psi(L_{\8})=\{p_0\}$, where $\psi$ is given by \eqref{eq:9b}, we deduce from \eqref{pecero} that $(h_x,h_y)=(\mu_{\rm max}, b_{\8})$, constant along $L_{\8}$. Also, observe that $(x_0,0)\in \cD\cap L_{\8}$ and recall that $\varphi(\cD)\subset \Pi_{\8}\cap \Sigma$. In this way, $\varphi(L_{\8})\subset \Pi_{\8}\cap \Sigma$. Since $h_x=\mu_{\rm max}$ along $L_{\8}$, we conclude by the definition of $\delta_0$ that $\delta_0\leq r_1$.

Consider next the geodesic arc $\beta^*$ in \eqref{eq:betastar}. It corresponds via \eqref{eq:9bb} to the half-line $\cL^*=\cL\cap \{y\geq \delta_0\}$. Since we have proved that $[r_1,r_2]\subset[\delta_0,\8)$, this geodesic arc $\beta^*$ has angle $\theta(\beta^*)$ greater than $\theta(\beta_{\8})= \theta_0$. This is a contradiction with the definition of $\theta_0$. Therefore, $(\nabla u)^{-1}(p_0)$ contains a geodesic semicircle of $\S^2$.

\vspace{0.2cm}

{\bf Step 9:} \emph{Existence of a geodesic semicircle in $(\nabla u)^{-1}(p)$ for at least $4$ different points.} 

We have seen that, for any Pogorelov point $p_0\in \nabla u(\S^2)$ of the hedgehog $\nabla u(\S^2)$, the set $(\nabla u)^{-1}(p_0)$ contains a geodesic semicircle. We will next show that there exist at least four different Pogorelov points for $\nabla u (\S^2)$, what proves the statement above.

Let $p$ be a contact point of $\nabla u(\S^2)$ with one of its support planes, and consider the set $\cN_p:=\{ \xi \in \S^2 : p \in \Pi_{\xi}\}$. Note that the convex hull $\cC$ of $\nabla u (\S^2)$ is not contained in a plane, since $\nabla u$ has some regular point of negative curvature (see the proof of Assertion \ref{ass:2}). In these conditions, it is well known that $\cN_p$ is a compact, convex subset of an open hemisphere of $\S^2$.

Arguing by contradiction, assume that $\nabla u(\S^2)$ has at most three (distinct) Pogorelov points $p_1,p_2,p_3$. Then $\cV:=\S^2\setminus \cup_{i=1}^3 \cN_{p_i}$ is a non-empty open set, since each $\cN_{p_i}$ lies in an open hemisphere. For almost any $\xi\in \cV$, the intersection $\Pi_{\xi}\cap \nabla u(\S^2)$ is a unique point $q_{\xi}$, which is not a Pogorelov point. Thus, from the definition of Pogorelov point, either $\nabla u (\xi)=q_{\xi}$, or $\nabla u (-\xi)=q_{\xi}$, for almost all $\xi \in \cV$. If for any such $\xi_0$ it holds $\nabla u(-\xi_0)\neq q_{\xi_0}$ , then, by definition of support plane, $$\esiz \nabla u(-\xi_0)-q_{\xi_0},\xi_0\esde <\esiz \nabla u(\xi_0)-q_{\xi_0},\xi_0\esde=0,$$ and so
 $$\esiz \nabla u(-\xi_0),\xi_0\esde <\esiz \nabla u(\xi_0),\xi_0\esde.$$ Hence, this property holds in a neighborhood $\cW\subset \cV$ of $\xi_0$, and it implies that for almost every $\xi\in \cW$, we have $\nabla u (\xi)=q_{\xi}$. In particular, $\nabla u$ is singular in a neighborhood of $\xi_0$, since regular points of $\nabla u(\S^2)$ never touch support planes. If $\nabla u(\xi_0)\neq q_{\xi_0}$, the same argument gives that $\nabla u$ is singular in a neighborhood $\cW$ of $-\xi_0$, and $\nabla u(\xi)=q_{-\xi}$ for almost every $\xi\in \cW$.

Finally, if $\nabla u (\xi)=\nabla u (-\xi)=q_{\xi}$ for almost all $\xi\in \cV$, we have that $\nabla u$ is singular in $\cV$. 

In other words, we have shown that there exists an open set $\cW\subset \S^2$ such that $\nabla u$ is singular everywhere on $\cW$, and for almost every $\xi\in \cW$, we have that $\nabla u (\xi)$ is the unique contact point of $\nabla u(\S^2)$ with one of the support planes $\Pi_{\xi}$ or $\Pi_{-\xi}$. 

Recall that, by homogeneity, $D^2 u$ always has a zero eigenvalue at every point, corresponding to the radial direction, and that the regular points of the hegdehog $\nabla u(\S^2)$ are those where the rank of $D^2 u$ is $2$; see the paragraph before Definition \ref{pogo}. Since $\nabla u$ is singular on $\cW$, by reducing $\cW$ if necessary, we can assume additionally that the rank of $D^2u$ is constantly equal to $0$ or $1$ in $\cW$. We rule out these two cases separately.

\begin{assertion}\label{ass:r2}
The rank of $D^2 u$ cannot be zero in $\cW$.
\end{assertion}
\begin{proof}
Assume that $D^2 u=0$ in $\cW$, and choose $\xi\in \cW$. Suppose, for definiteness, that $\nabla u (\xi)=q_{\xi}\in \Pi_{\xi}$; the discussion is similar if $\nabla u(\xi)\in \Pi_{-\xi}$. 

We will start arguing as in Step 3. Consider Euclidean coordinates $(x,y,z)$ in $\R^3$ such that $\xi=(1,0,0)$, and let $\Sigma$ be the entire saddle graph in $\R^3$ given by $z=h(x,y)$, where $h$ is defined by \eqref{eq:0}. Then, equations \eqref{ja3} and \eqref{ja3b} at the beginning of Step 3 hold, but \eqref{ja4} does not. Since $u$ is linear in a neighborhood of $\xi$, with $u_x=\mu_{\rm max}$, we deduce that instead of \eqref{ja4} we have in our context that
 \begin{equation}\label{ja4c}
 h_x(x,y) =\mu_{\rm max}  \hspace{1cm} \text{$\forall (x,y)\in (0,\8)\times \R$ with $x^2 \geq \delta (y^2+1),$}
 \end{equation}
for some $\delta>0$. In this way, if we choose $(x_0,0)$ with $x_0>\delta$ and define the linear function $$P(x,y):= \mu_{\rm max} (x-x_0) + h_y(x_0,0) y + h(x_0,0),$$  we have that $h(x,y)=P(x,y)$ in a connected planar subset $\Omega\subset \R^2$ that contains the set defined in \eqref{ja4c}, and $h(x,y)>P(x,y)$ in $\R^2-\Omega$.
 
 By the argument in Assertion \ref{lem4}, we deduce that $\R^2-\Omega$ is an open convex set. Consider the set $\Theta_0 \subset \S^2$ given by the points $\nu$ of the form \eqref{eq:9bb}, with $(x,y)\in \Omega$. Since \eqref{eq:9bb} is a totally geodesic mapping, this means that, if $\S_+^2:= \S^2\cap \{z>0\}$, then $\S_+^2 \setminus \Theta_0$ is a convex set of $\S_+^2$. But now, note that the Euclidean coordinates $(x,y,z)$ were chosen arbitrarily except for the condition $\xi=(1,0,0)$. Thus, if we define $\Theta\subset \S^2$ as the set of points $\nu\in \S^2$ that are given by \eqref{eq:9bb} for some $(x,y)\in \Omega$ with respect to \emph{some} Euclidean coordinates $(x,y,z)$ with $\xi=(1,0,0)$, we deduce then that $\S^2\setminus\Theta$ is a convex set of $\S^2$, and $u$ is linear on $\Theta$.  Then, $\S^2\setminus \Theta$ lies in an open hemisphere. Consequently, $u$ is linear on a closed hemisphere $H$ of $\S^2$, with $\nabla u=q_{\xi}$. Consider next the homogeneous function $v(p):=u(p)-\esiz p,q_{\xi}\esde$, defined for all $p\in \R^3$. Note that $D^2 v=D^2 u$ everywhere, and that $v$ vanishes along the geodesic $\parc H$ of $\S^2$. By \cite[Thm. 1.6.4]{NTV} or \cite[Thm. 2]{K}, $v$ must be linear. Hence, $u$ is linear, a contradiction.
 \end{proof}

\begin{assertion}\label{ass:r1}
The rank of $D^2 u$ cannot be $1$ in $\cW$.
\end{assertion}
\begin{proof}
In order to prove the assertion, we use some results of hegdehog theory developed by Martinez-Maure in \cite{MM0}, that we explain next. Given $h\in C^2(\S^2)$, let $\cH$ be the \emph{hedgehog} in $\R^3$ with support function $h$, i.e. $\cH$ is given by  $$\chi (\nu):= \nabla_{\S}\, h (\nu) + h(\nu) \nu : \S^2\flecha \cH\subset \R^3,$$ where $\nabla_{\S} $ denote the gradient in $\S^2$. We assume that the curvature of $\chi$ is negative at its regular points, and that $\chi$ is not constant. Note that the hedgehog $\cH:=\nabla u(\S^2)$ of our problem is in these conditions.

For any $\omega\in \S^2$, consider the plane $P=\{\omega\}^{\perp}$, and let $\pi:\R^3\flecha P$ denote the orthogonal projection. Define $\chi_{\omega}:  \S^1\equiv \S^2\cap P\flecha P$ by 
\begin{equation}\label{eripla}
\chi_{\omega}(\theta):= \pi(\chi (\theta)).
\end{equation}
Then, $\chi_{\omega}$ defines a \emph{planar hedgehog} in $P$, that we denote by $\cH_{\omega}=\chi_{\omega}(\S^1)$. Since $\cH$ has negative curvature at its regular points, this projected hedgehog $\cH_{\omega}$ has empty \emph{convex interior}; see Theorem 2 and Corollary 1 in \cite{MM0}, where the definition of convex interior of a planar hedgehog (which we will not use explicitly) is also presented; see also Corollary 1 in \cite{MM}. 

We now prove Assertion \ref{ass:r1} using this information. Since $D^2 u$ has rank one in the open set $\cW\subset \S^2$, then $\nabla u (\cW)$ is a regular curve $\gamma$. Also, note that for almost every $q\in\gamma$ we have either $\{q\}=\Pi_{\xi}\cap \nabla u(\S^2)$ or $\{q\}=\Pi_{-\xi}\cap \nabla u(\S^2)$.


%

Let $T$ be the unit tangent vector to $\gamma$ at $q$, and define $\omega:=T\times \xi$. Let $\pi:\R^3\flecha \{\omega\}^{\perp}$ denote the orthogonal projection onto $P=\{\omega\}^{\perp}$. Then $\beta:=\pi(\gamma)$ is a regular curve in $P=\{\omega\}^{\perp}$ around $\pi(q)$, and $\pi(q)\in \beta\cap \cH_{\omega}$ (since $\esiz T,\omega\esde=0$), where $\cH_{\omega}$ is the planar hedgehog given by \eqref{eripla}. Note that $\pi(q)$ is a regular point of $\cH_{\omega}$, since $\chi_{\omega}(T)=\pi(q)$ and $\esiz \nabla u(q), T\esde\neq 0$, by regularity of $\gamma$. Also, either $\cH_{\omega}$ lies on one side of the line $L_{\xi}=\Pi_{\xi}\cap P$, and in that case $\pi(q)\in L_{\xi}\cap \cH_{\omega}$, or else $\cH_{\omega}$ lies on one side of $L_{-\xi}=\Pi_{-\xi}\cap P$, and $\pi(q)\in L_{-\xi}\cap \cH_{\omega}$. In this way, in any of these two cases, the planar hedgehog $\cH_{\omega}\subset P$ touches one of its support lines at the regular point $\pi(q)$. Since $\cH_{\omega}$ has empty convex interior, we obtain a contradiction with \cite[Proposition 1]{MM0}. 
\end{proof}

Thus, we have proved that $\nabla u(\S^2)$ has at least four Pogorelov points, as claimed.

\vspace{0.2cm}

{\bf Step 10:} \emph{The final contradiction.} 

We now conclude the argument of the proof of Theorem \ref{th:main}. Recall that we had initially assumed that $u$ is not a linear function, and we were arguing by contradiction.

We have shown in Step 9 that there exist at least $4$ different points $p_1,\dots, p_4\in \nabla u(\S^2)$ for which $(\nabla u)^{-1}(p_j)$ contains a geodesic semicircle $\Gamma_j$ of $\S^2$. The geodesic semicircles $\Gamma_1,\dots, \Gamma_4$ are disjoint, since the $p_j$ are different.

Consider the region $\cO\subset \S^2$ defined below \eqref{coe}. By hypothesis on $\cO$, we have $\cO\cap \Gamma_j\neq \emptyset$ for some $j\in \{1,\dots, 4\}$. Let $\Omega_j$ denote the compact set $(\nabla u)^{-1}(p_j)$. Thus, $\Omega_j\cap \cO\neq \emptyset$ and, since $\cO$ is connected, either $\parc \Omega_j\cap \cO\neq \emptyset$ or $\cO\subset \Omega_j$.

Suppose, in the first place, that $\cO\subset \Omega_j$. Then, it is clear that the distance from $\cO$ to any of the semicircles $\Gamma_k$, $k\neq j$, is positive. In particular, there exists $\ep>0$ such that $\cO$ does not intersect the open set $\cU_{\ep}:=\{\nu \in \S^2: {\rm dist}(\nu,\Gamma_k)<\ep\}$. But on the other hand, it is clear that there exist infinitely many closed disjoint geodesic semicircles contained in $\cU_{\ep}$. This contradicts the hypothesis that $\cO$ intersects any configuration of $4$ disjoint geodesic semicircles. Thus, $\cO$ is not contained in $\Omega_j$.

Hence, there must exist some $w_j\in \parc \Omega_j\cap \cO$. Since $w_j\in (\nabla u)^{-1}(p_j)$, we can choose $w_j$ as the vector $\nu_0\in \S^2$ in the argument that we carried out in Steps 3 through 7. Specifically, choose Euclidean coordinates $(x,y,z)$ so that $\xi_j=(1,0,0)$ and $w_j=:\nu_0=(\nu_0^1,0,\nu_0^3)$, with $\nu_0^3>0$. Denote $\S^2_+=\S^2\cap \{z>0\}$. Then, by the argument in Steps 3 through 7, the connected component of the set $(\nabla u)^{-1}(p_j)\cap \S_+^2$ that contains $\nu_0$ is made of the points $\nu\in \S^2$ given by \eqref{eq:9bb}, with $(x,y)$ a point of the planar set $\cD$ defined in \eqref{eq:dee}. Also, \eqref{pecero} holds for $p_0:=p_j$.

Take $x_0\in \R$ given by $\nu_0=\frac{1}{\sqrt{1+x_0^2}}(x_0,0,1)$. Since $\nu_0\in \parc \Omega_j$, obviously $(x_0,0)\in \parc\cD$, and $h_x(x_0,0)=\mu_{\rm max}$ by \eqref{pecero} and \eqref{eq:9b}. Thus, $h_x$ has an absolute maximum at $(x_0,0)$. Hence, as $\nu_0:=w_j$ lies in $\cO$, it follows by Assertion \ref{ass:1} that $h_x$ is constant around $(x_0,0)$, since $\nabla h$ cannot be an open mapping. Then, by \eqref{eq:9b}, $\nu_0$ lies in the interior of $\Omega_j$, a contradiction with $\nu_0\in \parc \Omega_j$. 

By this final contradiction, the function $u$ must be linear, and this proves Theorem \ref{th:main}.
\section{Proof of Theorem \ref{th:4sh}}\label{sec:4s}

In Steps 2 through 9 of our proof of Theorem \ref{th:main} we actually showed the following result. Let $u\in C^2(\R^3\setminus \{0\})$ be a degree one homogeneous solution to a linear equation \eqref{edp1}. Assume that the coefficients $a_{ij}$ of \eqref{edp1} satisfy the degenerate ellipticity conditions (i), (ii) in \eqref{co}. Let $\nabla u:\S^2\flecha \R^3$ be the restriction of the gradient of $u$ to $\S^2$. Then, there exist at least $4$ different points $p_1,\dots, p_4$ in $\R^3$ such that each $(\nabla u)^{-1}(p_j)$ contains a geodesic semicircle $\Gamma_j$, for $j=1,\dots, 4$. These semicircles are disjoint, and $D^2u$ vanishes along the configuration $\Gamma=\cup_{i=1}^4 \Gamma_i$. 

As explained at the beginning of Step 2, there is an equivalence between degree one homogeneous solutions $u\in C^2(\R^3\setminus\{0\})$ of \eqref{edp1} whose coefficients satisfy conditions (i), (ii) in \eqref{co}, and $C^2$ saddle functions $v(x)=u(x/|x|)$ on $\S^2$. Taking into account this equivalence, it is then clear that the result obtained in Steps 2 through 9 that we just recalled directly proves Theorem \ref{th:4sh}.

Theorem \ref{th:4sh} is equivalent to the geometric statement below. Indeed, if $\rho\in C^2(\S^2)$ denotes the support function of an ovaloid satisfying \eqref{eqpan}, then $v:=\rho-c$ is a saddle function in $\S^2$, thus in the conditions of Theorem \ref{th:4sh} (and conversely).

\begin{theorem}\label{th:4s}
Let $S\subset \R^3$ be a $C^2$ ovaloid in $\R^3$ whose principal curvatures $\kappa_1,\kappa_2$ satisfy
 \begin{equation}\label{eqpan}
(\kappa_1-c)(\kappa_2-c)\leq 0
\end{equation} 
for some $c>0$. Then, $S$ is round along $4$ geodesic semicircles. Specifically, $S$ is tangent up to second order to four spheres $\Sigma_1^c,\dots, \Sigma_4^c$ of radius $1/c$ along four disjoint geodesic semicircles $\alfa_j\subset \Sigma_j^c\cap S$, for $j=1,\dots, 4.$
\end{theorem}
In other words, there exist $4$ disjoint geodesic semicircles $\Gamma_1,\dots, \Gamma_4$ in $\S^2$ such that, if $\eta:S\flecha \S^2$ is the Gauss map of $S$, then each $\eta^{-1}(\Gamma_j)=\alfa_j$ is made of umbilic points of $S$, and coincides with a geodesic semicircle of a sphere of radius $1/c$ in $\R^3$.

\def\refname{References}

\vskip 0.2cm

\noindent José A. Gálvez

\noindent Departamento de Geometría y Topología,\\ Instituto de Matemáticas IMAG, \\Universidad de Granada (Spain).

\noindent  e-mail: {\tt jagalvez@ugr.es}

\vskip 0.2cm

\noindent Pablo Mira

\noindent Departamento de Matemática Aplicada y Estadística,\\ Universidad Politécnica de Cartagena (Spain).

\noindent  e-mail: {\tt pablo.mira@upct.es}

\vskip 0.4cm

\noindent This research has been financially supported by: Projects PID2020-118137GB-I00 and CEX2020-001105-M, funded by MCIN/AEI /10.13039/501100011033, and Junta de Andalucia grants no. A-FQM-139-UGR18 and P18-FR-4049
\end{document}